# SECOND ORDER ASYMPTOTICS FOR MATRIX MODELS

By Alice Guionnet and Edouard Maurel-Segala

*Ecole Normale Supérieure de Lyon*

We study several-matrix models and show that when the potential is convex and a small perturbation of the Gaussian potential, the first order correction to the free energy can be expressed as a generating function for the enumeration of maps of genus one. In order to do that, we prove a central limit theorem for traces of words of the weakly interacting random matrices defined by these matrix models and show that the variance is a generating function for the number of planar maps with two vertices with prescribed colored edges.

**1. Introduction.** In this paper we study the asymptotics of Hermitian random matrices whose distribution is given by a small convex perturbation of the Gaussian Unitary Ensemble (denoted GUE). We shall consider $m$-tuples of random matrices, with an integer number $m \in \mathbb{N}$ fixed throughout this paper. Then, the law $\mu^N$ of $m$ independent matrices following the GUE is given, for $N \times N$ Hermitian matrices $\mathbf{A} = (A_1, \ldots, A_m)$, by

$$d\mu^N(\mathbf{A}) = e^{(-N/2)\text{Tr}(\sum_{i=1}^m A_i^2)} \prod_{i=1}^m \prod_{j=1}^N d(A_i)_{jj} \prod_{1 \le j < k \le N} d\Re e(A_i)_{jk} d\Im m(A_i)_{jk},$$

with Tr the nonnormalized trace $\text{Tr}(A) = \sum_{i=1}^N A_{ii}$. In other words, the $\mathbf{A} = (A_1, \ldots, A_m)$ are independent Hermitian matrices whose entries are, above the diagonal, independent complex centered Gaussian variables with variance $N^{-1}$. Let $V(\mathbf{X})$ be a polynomial in $m$ noncommutative indeterminates $\mathbf{X} = (X_1, \ldots, X_m)$ such that $\text{Tr}(V(\mathbf{A}))$ is real for all $m$-tuple of Hermitian matrices $\mathbf{A} = (A_1, \ldots, A_m)$. Then, we shall study the following probability measure $\mu_V^N$ on the set $\mathcal{H}_N(\mathbb{C})^m$ of $m$-tuple of $N \times N$ Hermitian matrices

$$d\mu_V^N(\mathbf{A}) = \frac{1}{Z_V^N} e^{-N\text{Tr}(V(\mathbf{A}))} d\mu^N(\mathbf{A}),$$

where $Z_V^N$ is the normalizing constant so that $\mu_V^N$ is a probability measure.

Besides, we require that the trace of $W(\mathbf{A}) := V(\mathbf{A}) + \frac{1}{2} \sum_{i=1}^m A_i^2$ is a strictly convex function of the entries of $\mathbf{A} = (A_1, \ldots, A_m) \in \mathcal{H}_N(\mathbb{C})^m$ for any $N \in \mathbb{N}$. In that case $Z_V^N$ is automatically finite. More precisely, for $c > 0$, we say that $V$ is $c$-convex if for any $N \in \mathbb{N}$, $\mathbf{A} \in \mathcal{H}_N(\mathbb{C})^m \to \mathrm{Tr}(W(\mathbf{A}))$ is real-valued and with Hessian bounded below by $cI$. An example of $c$-convex potential is

$$V(X_1, \ldots, X_m) = \sum_j P_j \left( \sum_i \alpha_i^j X_j \right) + \sum_{j,k} \beta_{j,k} X_i X_j,$$

with convex polynomials $P_j$ on $\mathbb{R}$, real numbers $\alpha_i^j, \beta_{j,k}$ and $\sum_j |\beta_{j,k}| \leq 1 - c$ for all $k \in \{1, \ldots, m\}$ (see Section 2 for more details).

The central result of this paper can roughly be stated as follows.

THEOREM 1.1.   *Let $V = V_\mathbf{t}(X_1, \ldots, X_m) = \sum_{j=1}^n t_j q_j(X_1, \ldots, X_m)$ be a polynomial potential with $n \in \mathbb{N}$, $\mathbf{t} = (t_1, \ldots, t_n) \in \mathbb{C}^n$ and monomials $(q_j)_{1 \leq j \leq n}$ fixed. For all $c > 0$, there exists $\eta > 0$ so that if $|\mathbf{t}| := \max_{1 \leq j \leq n} |t_j| \leq \eta$ and $V_\mathbf{t}$ is $c$-convex, there exists $F^i(V_\mathbf{t}) = F^i(t_1, \ldots, t_n)$ for $i = 0, 1$ so that*

$$\log Z_{V_\mathbf{t}}^N = N^2 F^0(V_\mathbf{t}) + F^1(V_\mathbf{t}) + o(1).$$

The first order expansion $F^0(V_\mathbf{t})$ was already obtained in [14] and we extend our study here to the second order. The higher order expansions can also be tackled by a refinement of our strategy; this is the subject of a separate article by Maurel-Segala [19]. Moreover, we believe our tools sufficiently robust to tackle other models such as the Gaussian orthogonal ensemble, or the Haar measure on the unitary group, for instance. Again, this is the subject of further studies.

We next turn to the combinatorial interpretation of $F^0(V_\mathbf{t})$, $F^1(V_\mathbf{t})$ has generating functions of maps.

Matrix models have been used intensively in physics in connection with the problem of enumerating maps; see the reviews [9, 12]. Let us recall that a map of genus $g$ is a graph which is embedded into a surface of genus $g$ in such a way that the edges do not intersect and dissecting the surface along the edges decomposes it into faces which are homeomorphic to a disk. We will call a star the couple of a vertex and the half-edges which are glued to this star. A star will have a distinguished half-edge and an orientation and will eventually have colored half-edges when $m \geq 2$. When $m = 1$, it is well known that if $V_\mathbf{t} = 0$ [i.e., $\mathbf{t} = (0, \ldots, 0)$] moments of the random matrices from the GUE are related with the enumeration of maps; for instance, the number $\mathcal{M}_k^g$ of maps with genus $g$ with one star with $2k$ half-edges were computed by Harer and Zagier [16] using the formula

$$\int \frac{1}{N} \text{Tr}(A_1^{2k}) \, d\mu^N(A_1) = \sum_{g=0}^{[k/2]} \frac{1}{N^{2k}} \mathcal{M}_k^g.$$

It was shown in [10] (see also [1, 2]) that when $m = 1$, this enumerative property extends to the free energy of matrix models at all orders, as conjectured and widely used in physics (see, e.g., [7]). More precisely, if $V_{\mathbf{t}} = \sum_{i=1}^{n} t_i x^{n_i}$ with $D = \max n_i = n_p$ even and $t_p / \sum_{i \neq p} |t_i|$ large enough, for all $k \in \mathbb{N}$, there exists $\eta > 0$ so that, for $|\mathbf{t}| \leq \eta$,

$$\log Z_{V_{\mathbf{t}}}^N = N^2 \sum_{g=0}^{k} \frac{1}{N^{2g}} F^g(V_{\mathbf{t}}) + o(N^{2-2k}),$$

with

$$F^g(V_{\mathbf{t}}) = \sum_{k_1, \ldots, k_n \in \mathbb{N}^n \setminus \{0, \ldots, 0\}} \prod_i \frac{(-t_i)^{k_i}}{k_i!} \mathcal{M}_{k_1, \ldots, k_n}^g,$$

where $\mathcal{M}_{k_1, \ldots, k_n}^g$ is the number of maps of genus $g$ with $k_i$ vertices of degree $n_i$, $1 \leq i \leq n$.

Several-matrices integrals are related with the enumeration of colored (or decorated) maps. To make this statement clear, let us associate to a monomial $q(X) = X_{i_1} \cdots X_{i_p}$ a colored star as follows. We choose $m$ different colors $\{1, \ldots, m\}$. The star associated to $q$ (called a star of type $q$) is a vertex equipped with colored half-edges such that the first half-edge has color $i_1$, the second has color $i_2$ till the last half-edge which has color $i_p$. Because the star has a distinguished half-edge (the one associated with $X_{i_1}$) and an orientation, this defines a bijection between noncommutative monomials and colored stars. Then, it can be seen [23] that, for any monomial $q$,

$$\lim_{N \to \infty} \int \frac{1}{N} \text{Tr}(q(\mathbf{A})) \, d\mu^N(\mathbf{A}) = \mathcal{M}_0(q),$$

with $\mathcal{M}_0(q)$ the number of planar maps with one-colored star of type $q$ such that only half-edges of the same color can be glued pair-wise together (then forming a one-colored edge). In [14], we proved that if $V_{\mathbf{t}}$ is $c$-convex and $\mathbf{t} = (t_1, \ldots, t_n)$ is small enough, the limit $F^0(V_{\mathbf{t}})$ of the free energy given in Theorem 1.1 is analytic in the variables $t_i$ in a neighborhood of the origin and its expansion is a generating function for planar maps with prescribed colored stars;

$$(1) \qquad F^0(V_{\mathbf{t}}) = \sum_{k_1, \ldots, k_n \in \mathbb{N}^n \setminus \{0, \ldots, 0\}} \prod_{i=1}^{n} \frac{(-t_i)^{k_i}}{k_i!} \mathcal{M}_{k_1, \ldots, k_n},$$

with $\mathcal{M}_{k_1,\ldots,k_n}$ the number of planar maps with $k_i$ colored stars of type $q_i$, the gluing being allowed only between half-edges of the same color. Note, however, that we cannot retrieve all the numbers $\mathcal{M}_{k_1,\ldots,k_n}$ from the $F^0(V_{\mathbf{t}})$'s because the condition that $\mathrm{Tr}(V_{\mathbf{t}})$ is real requires that the parameters $\mathbf{t}$ satisfy some relations. Namely, if $*$ denotes the involution $(zX_{i_1}\cdots X_{i_k})^* = \bar{z}X_{i_k}\cdots X_{i_1}$, we must have $\mathrm{Tr}(V_{\mathbf{t}}) = \frac{1}{2}\mathrm{Tr}(V_{\mathbf{t}} + V_{\mathbf{t}}^*)$ and, therefore, if $V_{\mathbf{t}} = \sum t_i q_i$, to each $t_i$ must corresponds a $t_j$ such that $\mathrm{Tr}(q_j) = \mathrm{Tr}(q_i^*)$ and $t_j = \bar{t_i}$. Thus, the $F^0(V_{\mathbf{t}})$'s are generating functions for the number of planar maps with $k_i$ colored stars of type $q_i$ or $q_i^*$. The convexity assumption also should induce some extra relations between the parameters, but it can be removed as shown in Theorem 1.4.

In this paper we shall prove that such a representation also holds for the correction $F^1(V_{\mathbf{t}})$ to the free energy given in Theorem 1.1.

PROPERTY 1.2.   $F^1(V_{\mathbf{t}})$ is analytic in the parameters $t_i$ in some neighborhood of the origin. Its expansion is a generating function of maps:

$$F^1(V_{\mathbf{t}}) = \sum_{k_1,\ldots,k_n \in \mathbb{N}^n \setminus \{0,\ldots,0\}} \prod_{i=1}^{n} \frac{(-t_i)^{k_i}}{k_i!}\mathcal{M}^1_{k_1,\ldots,k_n},$$

with $\mathcal{M}^1_{k_1,\ldots,k_n}$ the number of maps with genus one with $k_i$ colored stars of type $q_i$. In particular, the above sum converges absolutely for $\max_{1 \le i \le n}|t_i|$ small enough.

Let us remark that such a representation is commonly assumed to hold in physics since the formal result is always true for finite $N$. For a few models (viz., models similar to the Ising model on random graphs), the analysis has been pushed forward to actually give a rather explicit formula for the generating function $F^1(V_{\mathbf{t}})$ in terms of the limiting spectral measure of one matrix under the Gibbs measure $\mu^N_{V_{\mathbf{t}}}$ (see, e.g., Eynard, Kokotov and Korotkin [11]). Our strategy is here to study the most general potentials, providing a general formula for $F^1(V_{\mathbf{t}})$ in terms of the limiting empirical measure of all the matrices (see Section 6).

Our arguments to prove Theorem 1.1 are rather different from [10] or [1] where orthogonal polynomials were used. In [10], the idea was to develop a Riemann–Hilbert approach based on precise asymptotics of orthogonal polynomials. In the case of several-matrices models, the technology of orthogonal polynomials is far to be as much developed (except for the Ising model; see [6]). We shall therefore use different tools; the first, which is well spread in physics, is the use of the Schwinger–Dyson equation, the second, for which we need a convex potential, is the a priori concentration inequalities. To sketch our strategy, let us denote $\hat{\mu}^N$ the empirical measure

$$\hat{\mu}^N : P \longrightarrow \frac{1}{N}\mathrm{Tr}(P(\mathbf{A})) = \frac{1}{N}\mathrm{Tr}(P(A_1,\ldots,A_m)),$$

where $P$ runs over the set $\mathbb{C}\langle X_1, \ldots, X_m \rangle$ of noncommutative polynomials in $m$ indeterminates. Note that when $m = 1$, $\hat{\mu}^N$ is the spectral measure of $A_1$, and therefore a probability measure on $\mathbb{R}$. When $m \geq 2$, $\hat{\mu}^N$ is a tracial state, which generalizes the notion of measures to a noncommutative setting (see, e.g., [24]). Observe that, for $1 \leq i \leq m$,

$$\partial_{t_i} \log Z^N_{V_{\mathbf{t}}} = -N^2 \mu^N_{V_{\mathbf{t}}}(\hat{\mu}^N(q_i))$$

so that the second order asymptotics of the free energy will follow from that of $\overline{\mu}^N = \mu^N_{V_{\mathbf{t}}}[\hat{\mu}^N]$ evaluated at the monomials $q_i$, $1 \leq i \leq n$. Then, a simple integration by parts shows that, for any $N \in \mathbb{N}$, the following finite $N$ Schwinger–Dyson equation holds

$$\mu^N_{V_{\mathbf{t}}}(\hat{\mu}^N \otimes \hat{\mu}^N(\partial_i P)) = \mu^N_{V_{\mathbf{t}}}(\hat{\mu}^N((X_i + D_i V_{\mathbf{t}})P))$$

for any polynomial $P$ and $i \in \{1, \ldots, m\}$. Here, $\partial_i, D_i$ are noncommutative derivatives (see Section 2 for a definition). Based on this equation and concentration inequalities (see Section 2 for a definition). Based on this equation and concentration inequalities it was shown in [14] that for sufficiently small parameters $\mathbf{t} = (t_1, \ldots, t_n)$, $\hat{\mu}^N$ converges almost surely and in expectation (for the weak topology generated by the set $\mathbb{C}\langle X_1, \ldots, X_m \rangle$ of noncommutative polynomials). Its limit $\mu_{\mathbf{t}}$ is a solution of the Schwinger–Dyson equation

$$(2) \qquad \mu_{\mathbf{t}} \otimes \mu_{\mathbf{t}}(\partial_i P) = \mu_{\mathbf{t}}((X_i + D_i V_{\mathbf{t}})P)$$

$$\forall P \in \mathbb{C}\langle X_1, \ldots, X_m \rangle, 1 \leq i \leq m.$$

It is the unique solution which satisfies a bound of the form $|\mu_{\mathbf{t}}(X_i^d)| \leq C^d$ for all $d \in \mathbb{N}$ and all $i \in \{1, \ldots, m\}$, when $\mathbf{t} = (t_1, \ldots, t_n)$ is small enough and $C$ finite, independent of $\mathbf{t}$.

In this paper we investigate the correction to this convergence by proving a central limit theorem for $\hat{\mu}^N - \mu_{\mathbf{t}}$. More precisely if we define $\hat{\delta}^N_{\mathbf{t}}(P) := N(\hat{\mu}^N(P) - \mu_{\mathbf{t}}(P))$, then we show the following:

THEOREM 1.3.  *For all $c > 0$, there exists $\eta > 0$ such that for all $\mathbf{t}$ in $B_{\eta,c} = B(0, \eta) \cap \{\mathbf{t} | V_{\mathbf{t}} \text{ is } c\text{-convex}\}$, for all $P$ in $\mathbb{C}\langle X_1, \ldots, X_m \rangle$, under $\mu^N_{V_{\mathbf{t}}}$, $\hat{\delta}^N_{\mathbf{t}}(P)$ converges in law toward a complex centered Gaussian law $\gamma_P$. Moreover, $\{\gamma_P | P \in \mathbb{C}\langle X_1, \ldots, X_m \rangle\}$, equipped with the natural addition $\gamma_P + \gamma_Q = \gamma_{P+Q}$, is a Gaussian space and the covariance function is a generating function for planar maps with two prescribed stars.*

Such a central limit theorem was proved for more general potentials when $m = 1$ by Johansson in [18]. When $m = 1$ but the entries are not Gaussian, we refer the reader to [3]. In the case $m \geq 2$ but $V = 0$, the central limit theorem was obtained in [8, 13, 20]. Our proof is rather close to that of [18] and in the physics spirit; by doing an infinitesimal change of variables, it can be seen that the random variable

$$\hat{\delta}_{\mathbf{t}}^N(\xi_{\mathbf{t}}P) := \sum_{i=1}^m \hat{\delta}_{\mathbf{t}}^N((I \otimes \mu_{\mathbf{t}} + \mu_{\mathbf{t}} \otimes I)(\partial_i D_i P) - (X_i + D_i V_{\mathbf{t}})D_i P)$$

converges in law toward a centered Gaussian variable. The main issue is then to show that the $\xi_{\mathbf{t}}P$'s are dense in the set of polynomials. When $m = 1$, Johansson could use finite Hilbert transformation to invert the operator $\xi_{\mathbf{t}}$. In our case, we deal with a differential operator acting on noncommutative test functions and we prove by hand that it is invertible for sufficiently small $t_i$'s in Section 4. Clearly, our analysis is perturbative at this point and does not try to find the optimal domain of validity of the central limit theorem.

To use the central limit theorem to obtain the second order asymptotics of $\mu_{V_{\mathbf{t}}}^N(\hat{\delta}_{\mathbf{t}}^N(P))$ observe that, by the finite dimensional Schwinger–Dyson equation, we get

$$N\mu_{V_{\mathbf{t}}}^N(\hat{\delta}_{\mathbf{t}}^N(\xi_{\mathbf{t}}P)) = \mu_{V_{\mathbf{t}}}^N(\hat{\delta}_{\mathbf{t}}^N \otimes \hat{\delta}_{\mathbf{t}}^N(\partial_i D_i P))$$

and the right-hand side converges toward the variance of the central limit theorem. So again, to obtain the limit of $N\mu_{V_{\mathbf{t}}}^N(\hat{\delta}_{\mathbf{t}}^N(P))$, we need to invert the operator $\xi_{\mathbf{t}}$ (see Section 6). The resulting formula for the free energy and the variance are given in terms of differential operators acting on noncommutative polynomial functions. Note that a similar formula for the variance of the central limit theorem governing the fluctuations of words of band matrices was found in [13]. Their interpretation in terms of enumeration of maps can be retrieved from the interpretation of noncommutative derivatives in terms of natural operations on maps (see [14]).

Finally, in the spirit of [14], we study matrix models with a nonnecessarily convex potential $V$. Since in that case $Z_V^N$ has no reason to be finite, we need to add a cut-off. For a positive constant $L$, we define

$$\mu_{V,L}^N(d\mathbf{A}) = \frac{1}{Z_{V,L}^N} \mathbb{1}_{\lambda_{\max}(\mathbf{A}) < L} \, e^{-N\mathrm{Tr}(V(\mathbf{A}))} \, d\mu^N(\mathbf{A}),$$

with $\lambda_{\max}(\mathbf{A})$ the maximum of the spectral radius of the $A_i$'s and $Z_{V,L}^N$ a normalizing constant. The remarkable point that we shall prove is that asymptotically the behavior of this measure is independent of $L$ and gives the same type of expansion as in the convex case.

THEOREM 1.4. Let $V = V_{\mathbf{t}}(X_1, \ldots, X_m) = \sum_{j=1}^n t_j q_j(X_1, \ldots, X_m)$ be a polynomial potential with $n \in \mathbb{N}$, $\mathbf{t} = (t_1, \ldots, t_n) \in \mathbb{C}^n$ and monomials $(q_j)_{1 \le j \le n}$ fixed. Assume that $\mathrm{Tr}(V_{\mathbf{t}}(\mathbf{A})))$ is real for all $\mathbf{A} \in \mathcal{H}_N(\mathbb{C})^m$, all $N \in \mathbb{N}$. There exists $L_0 > 0$ such that, for all $L > L_0$, there exists $\eta > 0$ so that if $|\mathbf{t}| := \max_{1 \le j \le n} |t_j| \le \eta$ then

$$\log Z_{V_{\mathbf{t}},L}^N = N^2 F^0(V_{\mathbf{t}}) + F^1(V_{\mathbf{t}}) + o(1)$$

with $F^0(V_{\mathbf{t}}), F^1(V_{\mathbf{t}})$ *as in* *(1)* *and Property* 1.2.

Note that in the large $N$ limit, the dependence in $L$ disappears.

In the next section we will describe our hypothesis of convexity and show some useful consequences. In Section 3 we give an estimate on the rate of convergence of $\mu_{V_{\mathbf{t}}}^N[\hat{\mu}^N]$ to $\mu_{\mathbf{t}}$. Then, in Section 4 we prove a central limit theorem, first only for some specific polynomials and then for arbitrary polynomials. In Sections 5 and 6, we give an interpretation of the variance and of the free energy in terms of enumeration of maps. Finally, in Section 7 we give some hints to generalize our proofs to the setting of Theorem 1.4.

## 2. Convex hypothesis and standard consequences.

### 2.1. *Framework and standard notation.*

2.1.1. *Noncommutative polynomials.* We denote $\mathbb{C}\langle X_1, \ldots, X_m \rangle$ the set of complex polynomials on the noncommutative unknown $X_1, \ldots, X_m$. Let $*$ denote the linear involution such that for all complex $z$ and all monomials

$$(zX_{i_1} \cdots X_{i_p})^* = \overline{z} X_{i_p} \cdots X_{i_1}.$$

We will say that a polynomial $P$ is self-adjoint if $P = P^*$ and denote $\mathbb{C}\langle X_1, \ldots, X_m \rangle_{sa}$ the set of self-adjoint elements of $\mathbb{C}\langle X_1, \ldots, X_m \rangle$.

For an integer number $N$, we denote $\mathcal{H}_N(\mathbb{C})$ the set of $N \times N$ Hermitian matrices. We shall sometimes identify $\mathcal{H}_N(\mathbb{C})$ with the set $\mathbb{R}^{N^2}$ of the corresponding real entries [by the bijection which associates to $A \in \mathcal{H}_{N(\mathbb{C})}$ the $N^2$-tuple $((\Re e(A_{ij})_{1 \le i \le j \le N}, (\Im m(A_{ij})_{1 \le i < j \le N}))]$.

Moreover, we shall denote in general by $A$ a random matrix, by $X$ a generic noncommutative indeterminate (e.g., to write polynomials). Bold symbols will in general denotes vectors; $\mathbf{A}$ (resp. $\mathbf{X}$) will in general denote a $m$-tuple of matrices (resp. noncommutative indeterminates), whereas $\mathbf{t}$ will denote a vector of complex scalars.

The potential $V$ will be later on assumed to be self-adjoint which guarantees that, for all integer $N$, all $\mathbf{A} = (A_1, \ldots, A_m) \in \mathcal{H}_N(\mathbb{C})^m$, $\mathrm{Tr}(V(\mathbf{A}))$ is real. Note that, conversely, if $\mathrm{Tr}(V(\mathbf{A}))$ is real, $\mathrm{Tr}(V(\mathbf{A})) = \mathrm{Tr}((V + V^*)(\mathbf{A})/2)$ and so we can replace $V$ by $(V + V^*)/2$ without loss of generality.

We shall assume also that $V$ satisfies some convexity property in this paper. Namely, we will say that $V$ is convex if, for any $N \in \mathbb{N}$,

$$\phi_V^N : \mathcal{H}_N(\mathbb{C})^m \simeq (\mathbb{R}^{N^2})^m \longrightarrow \mathbb{R},$$

$$((A_k)_{ij})_{\substack{1 \le i \le j \le N \\ 1 \le k \le m}} \longrightarrow \mathrm{Tr}(V(A_1, \ldots, A_m))$$

is a convex function of its entries.

Note that as we add a Gaussian potential $\frac{1}{2}\sum_{i=1}^m X_i^2$ to $V$ we can relax the hypothesis a little. We will say that $V$ is $c$-convex if $c > 0$ and $V + \frac{1-c}{2}\sum_1^m X_i^2$ is convex. Then the Hessian of $\phi_W^N$ with $W = V + \frac{1}{2}\sum_1^m X_i^2$ is symmetric positive with eigenvalues bigger than $c$.

An example is

$$V_{\mathbf{t}}(X_1, \ldots, X_m) = \sum_{i=1}^n P_i\left(\sum_{k=1}^m \alpha_k^i X_k\right) + \sum_{k,l} \beta_{k,l} X_k X_l,$$

with convex real polynomials $P_i$ in one unknown and real $\alpha_k^i, \beta_{k,l}$ such that, for all $l$, $\sum_k |\beta_{k,l}| \leq (1-c)$. This is due to Klein's lemma (see [15]) which states that the trace of a real convex function of a self-adjoint matrix is a convex function of the entries of the matrix.

In the rest of the paper we shall assume that $V$ is $c$-convex for some $c > 0$ fixed. We will denote $B(0, \eta) = \{\mathbf{t} \in \mathbb{C}^n : \max_{1 \leq i \leq n} |t_i| \leq \eta\}$ and $B_{\eta,c} = B(0,\eta) \cap \{\mathbf{t} : V_{\mathbf{t}} \text{ is } c\text{-convex}\}$.

2.1.2. *Noncommutative derivatives.* We define for $1 \leq i \leq m$ the noncommutative derivatives $\partial_i : \mathbb{C}\langle X_1, \ldots, X_m \rangle \to \mathbb{C}\langle X_1, \ldots, X_m \rangle^{\otimes 2}$ by the Leibniz rule

$$\partial_i PQ = \partial_i P \times (1 \otimes Q) + (P \otimes 1) \times \partial_i Q$$

and $\partial_i X_j = \mathbb{1}_{i=j} 1 \otimes 1$. So for a monomial $P$, the following holds:

$$\partial_i P = \sum_{P = RX_i S} R \otimes S,$$

where the sum runs over all possible monomials $R, S$ so that $P$ decomposes into $RX_i S$. We can iterate the noncommutative derivatives; for instance, $\partial_i^2 : \mathbb{C}\langle X_1, \ldots, X_m \rangle \to \mathbb{C}\langle X_1, \ldots, X_m \rangle \otimes \mathbb{C}\langle X_1, \ldots, X_m \rangle \otimes \mathbb{C}\langle X_1, \ldots, X_m \rangle$ is given on monomial functions by

$$\partial_i^2 P = 2 \sum_{P = RX_i SX_i Q} R \otimes S \otimes Q.$$

We denote $\sharp : \mathbb{C}\langle X_1, \ldots, X_m \rangle^{\otimes 2} \times \mathbb{C}\langle X_1, \ldots, X_m \rangle \to \mathbb{C}\langle X_1, \ldots, X_m \rangle$ the map $P \otimes Q \sharp R = PRQ$ and generalize this notation to $P \otimes Q \otimes R \sharp(S, T) = PSQTR$. So $\partial_i P \sharp R$ corresponds to the derivative of $P$ with respect to $X_i$ in the direction $R$, and, similarly, $2^{-1}[D_i^2 P \sharp(R, S) + D_i^2 P \sharp(S, R)]$ the second derivative of $P$ with respect to $X_i$ in the directions $R, S$.

We also define the so-called cyclic derivative $D_i$. If $m$ is the map $m(A \otimes B) = BA$, let us define $D_i = m \circ \partial_i$. For a monomial $P$, $D_i P$ can be expressed as

$$D_i P = \sum_{P = RX_i S} SR.$$

We shall denote in short $\mathbf{D}$ the cyclic gradient $(D_1, \ldots, D_m)$.

2.1.3. *Noncommutative laws.* For $(A_1, \ldots, A_m) \in \mathcal{H}_N(\mathbb{C})^m$, we define the linear form $\hat{\mu}^N_{A_1, \ldots, A_m}$ on $\mathbb{C}\langle X_1, \ldots, X_m \rangle$ by

$$\hat{\mu}^N_{A_1, \ldots, A_m}(P) = \frac{1}{N} \mathrm{Tr}(P(A_1, \ldots, A_m)),$$

with Tr the standard trace $\mathrm{Tr}(A) = \sum_{i=1}^N A_{ii}$. $\hat{\mu}^N_{A_1, \ldots, A_m}$ will sometimes be called the empirical distribution of the matrices $(A_1, \ldots, A_m)$. When there is no ambiguity and the matrices $A_1, \ldots, A_m$ follow the law $\mu^N_V$, we shall drop the subscript $A_1, \ldots, A_m$; $\hat{\mu}^N = \hat{\mu}^N_{A_1, \ldots, A_m}$. In [14] it was shown that if $V_{\mathbf{t}} = \sum_i t_i q_i$ is $c$-convex, for $|\mathbf{t}| := \max_{1 \le i \le n} |t_i|$ small enough, $\hat{\mu}^N$ converges weakly in expectation and almost surely under $\mu^N_V$ toward a limit $\mu_{\mathbf{t}}$ [i.e., for all $P$ in $\mathbb{C}\langle X_1, \ldots, X_m \rangle$, $\hat{\mu}^N(P)$ converges in expectation and almost surely to $\mu_{\mathbf{t}}(P)$]. We denote

$$\overline{\mu}^N_{\mathbf{t}}(P) = \mu^N_{V_{\mathbf{t}}}[\hat{\mu}^N(P)].$$

We shall later estimate differences between $\hat{\mu}^N$ and its limit. So, we set

$$\hat{\delta}^N_{\mathbf{t}} = N(\hat{\mu}^N - \mu_{\mathbf{t}}),$$

$$\overline{\delta}^N_{\mathbf{t}} = \int \hat{\delta}^N \, d\mu^N_V = N(\overline{\mu}^N_{\mathbf{t}} - \mu_{\mathbf{t}}),$$

$$\hat{\underline{\delta}}^N_{\mathbf{t}} = N(\hat{\mu}^N - \overline{\mu}^N_{\mathbf{t}}) = \hat{\delta}^N_{\mathbf{t}} - \overline{\delta}^N_{\mathbf{t}}.$$

In order to simplify the notation, we will make $\mathbf{t}$ implicit and drop the subscript $\mathbf{t}$ in the rest of this paper so that we will denote $\overline{\mu}^N, \mu, \hat{\delta}^N, \overline{\delta}^N$ and $\hat{\underline{\delta}}^N$ in place of $\overline{\mu}^N_{\mathbf{t}}, \mu_{\mathbf{t}}, \hat{\delta}^N_{\mathbf{t}}, \overline{\delta}^N_{\mathbf{t}}$ and $\hat{\underline{\delta}}^N_{\mathbf{t}}$, as well as $V$ in place of $V_{\mathbf{t}}$.

2.2. *Brascamp–Lieb inequality and a priori controls.* We use here a generalization of the Brascamp–Lieb inequality shown by Hargé [17] which implies that if $V$ is $c$-convex, for all convex function $g$ on $(\mathbb{R})^{mN^2} \simeq \mathcal{H}_N(\mathbb{C})^m$,

$$(3) \qquad \int g(\mathbf{A} - \mathbf{M}) \, d\mu^N_V(\mathbf{A}) \le \int g(\mathbf{A}) \, d\mu^N_c(\mathbf{A}),$$

where $\mathbf{M} = \int \mathbf{A} \, d\mu^N_V(\mathbf{A})$ is the $m$-tuple of deterministic matrices with entries $(\mathbf{M}_j)_{k\ell} = \int (A_j)_{k\ell} \, d\mu^N_V(\mathbf{A})$ for $k, \ell \in \{1, \ldots, N\}$, $j \in \{1, \ldots, m\}$, and $\mu^N_c$ is the Gaussian law on $\mathcal{H}_N(\mathbb{C})^m$ with covariance $(Nc)^{-1}$, that is, $\int f(\mathbf{A}) \, d\mu^N_c(\mathbf{A}) = \int f(c^{-1/2}\mathbf{A}) \, d\mu^N(\mathbf{A})$ for all measurable function $f$ on $\mathbb{R}^{mN^2}$.

Recall that $B_{\eta,c}$ is the subset of the complex numbers $\mathbf{t} \in \mathbb{C}^n$ which are bounded by $\eta$ and so that $V$ is $c$-convex. Based on the Brascamp–Lieb inequality, the following was shown in [14] (Theorem 3.4):

LEMMA 2.1 (Compact support). *If $c, \eta > 0$, then there exists $C_0 = C_0(c, \eta)$ finite such that, for all $i \in \{1, \ldots, m\}$, all $n \in \mathbb{N}$, all $\mathbf{t} \in B_{\eta,c}$,*

$$\mu(X_i^{2n}) \le \limsup_N \overline{\mu}^N(X_i^{2n}) \le C_0^{2n}.$$

Note that this lemma shows that, for $i \in \{1, \ldots, m\}$, the spectral measure of $X_i$ is asymptotically contained in the compact set $[-C_0, C_0]$.

PROOF OF LEMMA 2.1. Let us recall the proof of this result for completeness. Let $k \in \{1, \ldots, m\}$. As $g \colon \mathbf{A} \in \mathcal{H}_N(\mathbb{C})^m \to N^{-1} \mathrm{Tr}(A_k^{4d}) = \hat{\mu}^N(X_k^{4d})$ is convex by Klein's lemma, we can use the Brascamp–Lieb inequalities (3) to see that

$$\overline{\mu}^N((X_k - M_k)^{4d}) \leq \mu_c^N(\hat{\mu}^N(X_k^{4d})), \tag{4}$$

where $M_k := \mu_V^N(A_k)$ is the deterministic matrix with entries $(M_k)_{ij} = \mu_V^N \times ((A_k)_{ij})$. Thus, since $\mu_c^N(\hat{\mu}^N(X_i^{4d}))$ converges by Wigner theorem [25] toward $c^{-2d} C_{2d} \leq (c^{-1}4)^{2d}$ with $C_{2d}$ the Catalan number, we only need to control $M_k$. First observe that, for all $k$, the law of $A_k$ is invariant under the unitary group so that, for all unitary matrices $U$,

$$M_k = \mu_V^N[U A_k U^*] = U \mu_V^N[A_k] U^* \Rightarrow M_k = \mu_V^N(\hat{\mu}^N(X_k))I = \overline{\mu}^N(X_k)I. \tag{5}$$

Let us bound $\overline{\mu}^N(X_k)$. Jensen's inequality implies

$$Z_N^V \geq e^{-N^2 \mu^N((1/N)\mathrm{Tr}(V))} = e^{-N^2 \mu^N \circ \hat{\mu}^N(V)}.$$

According to [23], $\mu^N \circ \hat{\mu}^N$ converges in moments to the law of $m$ free semicircular operators, which are uniformly bounded. Thus, there exists a finite constant $L$ such that $Z_N^V \geq e^{-N^2 L}$. We now use the convexity of $V$ to find that, for all $N$, all $\mathbf{A} = (A_1, \ldots, A_m) \in \mathcal{H}_N(\mathbb{C})^m$,

$$\mathrm{Tr}\left(V(\mathbf{A}) + \frac{1-c}{2}\sum_{i=1}^m A_i^2\right) \geq \mathrm{Tr}\left(V(0) + \sum_{i=1}^m D_i V(0) A_i + (1-c)\sum_{i=1}^m A_i\right).$$

By Chebyshev's exponential inequality, and then using the above bound, we therefore obtain that, for any $\lambda \geq 0$,

$$\begin{aligned}
\mu_V^N(\hat{\mu}^N(X_k) \geq y) &\leq e^{-\lambda N^2 y} \mu_V^N(e^{\lambda N^2 \hat{\mu}^N(X_k)}) \\
&= e^{-\lambda N^2 y} \frac{Z_c^N}{Z_1^N Z_V^N} \mu_c^N(e^{\lambda N^2 \hat{\mu}^N(X_k) - N\mathrm{Tr}(V(\mathbf{A}) + ((1-c)/2)\sum_{i=1}^m A_i^2)}) \\
&\leq e^{N^2(L - V(0) - \lambda y + m/2 \log c)} \mu_c^N \\
&\quad \times (e^{-N\mathrm{Tr}(\sum_{i=1}^m ((1-c) + D_i V(0))A_i - \lambda A_k)}) \\
&= e^{N^2(L - V(0) - \lambda y + m/2 \log c)} \\
&\quad \times e^{(N^2/(2c))\sum_{i \neq k}(1-c+D_i V(0))^2 + (N^2/(2c))(1-c+D_k V(0)-\lambda)^2},
\end{aligned}$$

where we denoted $V(0) := V(0, \ldots, 0)$ and $D_i V(0) := D_i V(0, \ldots, 0)$. Remark that these constants are uniformly bounded for $\mathbf{t}$ in $B(0, R)$, $R > 0$. Thus, we deduce that

$$\mu_V^N(\hat{\mu}^N(X_k) \geq y) \leq e^{N^2[(a - \lambda y) + (1/(2c))(\lambda - b)^2]}$$

with two constants $a, b$ which are uniformly bounded in terms of $c, \eta$ for $\mathbf{t} \in B_{\eta,c}$. Optimizing with respect to $\lambda$ shows that there exists $A < +\infty$ so that, for $\mathbf{t}$ in $B_{\eta,c}$,

$$\mu_V^N(\hat{\mu}^N(X_k) \geq y) \leq e^{N^2(a-(c/2)y^2-by)}$$
$$\leq e^{N^2(A-(c/4)y^2)}.$$

Replacing $X_k$ by $-X_k$, we bound similarly $\mu_V^N(\hat{\mu}^N(X_k) \leq -y)$ and, hence, we have proved

$$\mu_V^N(|\hat{\mu}^N(X_k)| \geq y) \leq 2e^{N^2(A-(c/4)y^2)}.$$

As a consequence,

$$\mu_V^N(|\hat{\mu}^N(X_k)|) = \int_0^\infty \mu_V^N(|\hat{\mu}^N(X_k)| \geq y)\, dy$$
$$\leq 2\sqrt{c^{-1}A} + 2\int_{2\sqrt{c^{-1}A}}^\infty e^{-(N^2c/4)(y^2-4A/c)}\, dy \leq 4\sqrt{c^{-1}A}, \tag{6}$$

where the last inequality holds for $N$ sufficiently large. Recall that $A$ is a continuous function of the $t_i$'s and, therefore, our bound on $\sup_N \mu_V^N(|\hat{\mu}^N(X_k)|)$, which controls the spectral radius of $M_k$ in any dimension $N$, is locally bounded in $\mathbf{t}$. This completes the proof with (4). □

Let us derive some other useful properties due to the convexity hypothesis. Let $\lambda_{\max}^N(A_i)$ be the maximum of the absolute value of the eigenvalues of $A_i$. We first obtain an estimate on $\lambda_{\max}^N(\mathbf{A})$, the maximum of the $(\lambda_{\max}^N(A_i))_{1 \leq i \leq m}$ under the law $\mu_V^N$.

LEMMA 2.2 (Exponential tail of the largest eigenvalue). *If $c, \eta > 0$, then there exists $\alpha > 0$ and $M_0 < \infty$ such that, for all $\mathbf{t} \in B_{\eta,c}$, all $M \geq M_0$ and all integer $N$,*

$$\mu_V^N(\lambda_{\max}^N(\mathbf{A}) > M) \leq e^{-\alpha MN}.$$

PROOF. Since the largest eigenvalue

$$\lambda_{\max}^N(\mathbf{A}) = \max_{1 \leq i \leq m} \sup_{\|u\|=1} \langle u, A_i A_i^* u \rangle^{1/2}$$

is a convex function of the entries of the $A_i$'s, we can apply the Brascamp–Lieb inequality (3) to obtain that, for all $s \in [0, \frac{c}{10}]$,

$$\int e^{sN\lambda_{\max}^N(\mathbf{A}-\mathbf{M})}\, d\mu_V^N(\mathbf{A}) \leq \int e^{sN\lambda_{\max}^N(\mathbf{A})}\, d\mu_c^N(\mathbf{A}) \leq C_0^N,$$

where the last inequality comes from the bound on the largest eigenvalue of the GUE shown, for instance, in [5]. Now,

$$\lambda_{\max}^N(\mathbf{A}) \le \lambda_{\max}^N(\mathbf{A} - \mathbf{M}) + \lambda_{\max}^N(\mathbf{M}) \le \lambda_{\max}^N(\mathbf{A} - \mathbf{M}) + 4\sqrt{Ac^{-1}},$$

where we used the bound (6). Consequently, we deduce that

$$\int e^{sN\lambda_{\max}^N(\mathbf{A})} \, d\mu_V^N(\mathbf{A}) \le C^N$$

for a positive finite constant $C$. We conclude by a simple application of Chebyshev's inequality.   $\square$

2.3. *Concentration inequalities.* We next turn to concentration inequalities for the trace of polynomials on the subset of $\mathcal{H}_N(\mathbb{C})^m \simeq \mathbb{R}^{N^2 m}$:

$$\Lambda_M^N = \left\{ \mathbf{A} \in \mathcal{H}_N(\mathbb{C})^m : \lambda_{\max}^N(\mathbf{A}) = \max_i(\lambda_{\max}^N(A_i)) \le M \right\}$$

for some fixed $M > 0$. Recall that $\underline{\hat{\delta}}^N = N(\hat{\mu}^N - \bar{\mu}^N)$. We shall prove concentration inequalities for $\underline{\hat{\delta}}^N(P)$ on $\Lambda_M^N$ for polynomial functions $P$. However, concentration inequalities should not restrict to polynomial functions but hold more generally for Lipschitz functions (see, e.g., [15]). We thus define the following Lipschitz semi-norm:

$$(7) \qquad \|P\|_{\mathcal{L}}^M = \sup_{\mathcal{A}C*\text{-algebra}} \sup_{\substack{x_1,\dots,x_m \in \mathcal{A} \\ \forall i, x_i = x_i^*, \|x_i\|_{\mathcal{A}} \le M}} \left( \sum_{k=1}^m \|D_k P D_k P^*\|_{\mathcal{A}} \right)^{1/2}.$$

Be aware that this is not a norm, since, for example, $\|1\|_{\mathcal{L}}^M = 0$ or $\|X_1 X_2 - X_2 X_1\|_{\mathcal{L}}^M = 0$. However, on these particular polynomials, $\underline{\hat{\delta}}^N$ vanishes. This fact can be generalized as follows; if we set

$$m_{M,P}^N := \frac{1}{\mu_V^N(\Lambda_M^N)} \mu_V^N(\mathbb{1}_{\Lambda_M^N} \underline{\hat{\delta}}^N(P)),$$

we shall see (see the proof below) that on $\Lambda_M^N$

$$|\underline{\hat{\delta}}^N(P) - m_{M,P}^N| \le 2M\sqrt{m}N\|P\|_{\mathcal{L}}^M.$$

Therefore, if we denote $\mathbb{C}\langle X_1,\dots,X_m\rangle_{\mathcal{L}}^M$ the completion and separation of $\mathbb{C}\langle X_1,\dots,X_m\rangle$ for $\|\cdot\|_{\mathcal{L}}^M$, we can extend $\underline{\hat{\delta}}^N - m_{M,P}^N$ to $\mathbb{C}\langle X_1,\dots,X_m\rangle_{\mathcal{L}}^M$ on $\Lambda_M^N$. A similar result will be proved for $\mu$ in Lemma 4.9 (note, however, that the arguments of this lemma do not apply here because $\bar{\mu}^N$ is not the law of uniformly bounded matrices).

We shall prove the following:

LEMMA 2.3 (Concentration inequality). *Let* **t** *be such that* $V$ *is* $c$-*convex. There exists* $\alpha, M_0 > 0$ *such that, for all* $N$ *in* $\mathbb{N}$, *all* $M > M_0$, *all* $P \in \mathbb{C}\langle X_1, \ldots, X_m \rangle_{\mathcal{L}}^M$, *there exists a positive constant* $\varepsilon_{P,M}^N$ *such that, for any* $\varepsilon > 0$,

$$(8) \qquad \mu_V^N(\{|\hat{\hat{\varrho}}^N(P) - m_{P,M}^N| \geq \varepsilon + \varepsilon_{P,M}^N\} \cap \Lambda_M^N) \leq 2e^{-c\varepsilon^2/(2(\|P\|_{\mathcal{L}}^M)^2)}.$$

*Moreover, there exists a universal constant* $C$ *such that*

$$\varepsilon_{P,M}^N \leq 2CNM\|P\|_{\mathcal{L}}^M e^{-\alpha NM}.$$

*If* $P$ *is a monomial of degree* $0 < d < \alpha N$, *we have*

$$\|P\|_{\mathcal{L}}^M \leq dM^{d-1}, \qquad \varepsilon_{P,M}^N \leq NCdM^d e^{-\alpha MN},$$

$$|m_{P,M}^N| \leq N(3M^d + d^2)e^{-\alpha MN}.$$

PROOF. Since $V$ is $c$-convex, for all integer number $N$, the Hessian of

$$\phi_V^N : \mathbf{A} \simeq ((A_k)_{ij})_{\substack{1 \leq i \leq j \leq N \\ 1 \leq k \leq m}} \in \mathbb{R}^{mN^2} \longrightarrow \mathrm{Tr}(V(A_1, \ldots, A_m)) \in \mathbb{R}$$

is bounded below by $cI$. Therefore, since $\mu_V^N$ has density $e^{-N\phi_V^N(\mathbf{A})}$ with respect to the Lebesgue measure, $\mu_V^N$ satisfies a Log–Sobolev inequality with constant $(Nc)^{-1}$ (see, e.g., Corollaire 5.5.2, page 87 in [4]). In other words, for any continuously differentiable function $f$ from $\mathbb{R}^{mN^2}$ into $\mathbb{R}$,

$$\int f^2 \log \frac{f^2}{\mu_V^N(f^2)} d\mu_V^N \leq \frac{2}{Nc} \int \|\nabla f\|^2 d\mu_V^N,$$

with $\nabla f$ the gradient of $f$ and $\|\cdot\|$ the Euclidean norm. Here and in the sequel we identify the measure $\mu_V^N$ as a measure on $\mathbb{R}^{N^2m}$. This implies, by the well-known Herbst argument (see, e.g., [4], Théorème 7.4.1, page 123), that $\mu_V^N$ satisfies concentration inequalities. Let us briefly summarize this argument for completeness. If $f$ is a continuously differentiable function, differentiating $X(\lambda) := \frac{1}{\lambda} \log \mu_V^N[e^{\lambda f}]$ and using the Log–Sobolev inequality yields

$$\partial_\lambda X(\lambda) \leq \frac{2}{cN\lambda^2 \mu_V^N(e^{\lambda f})} \mu_V^N(\|\nabla e^{(1/2)\lambda f}\|^2) \leq \frac{1}{2cN} \|\|\nabla f\|^2\|_\infty.$$

If we assume $\mu_V^N(f) = 0$, we find that $X(0) = 0$ and so integrating with respect to $\lambda$ yields

$$\mu_V^N(e^{\lambda f}) \leq e^{\lambda^2 \|\|\nabla f\|^2\|_\infty/(2cN)}.$$

Using Chebyshev's inequality thus gives, for $\varepsilon > 0$ and $\lambda > 0$,

$$\mu_V^N(f \geq \varepsilon) \leq e^{-\lambda\varepsilon} e^{\lambda^2 \|\|\nabla f\|^2\|_\infty/(2cN)}$$

and so optimizing with respect to $\lambda$ results with

$$\mu_V^N(f \geq \varepsilon) \leq e^{-cN\varepsilon^2/(2\|\|\nabla f\|\|_\infty^2)}.$$

Replacing $f$ by $-f$ gives the well-known concentration estimate, for any $\varepsilon > 0$,

$$\mu_V^N(|f| \geq \varepsilon) \leq 2e^{-cN\varepsilon^2/(2\|\|\nabla f\|\|_\infty^2)}$$

for any continuously differentiable function $f$ such that $\mu_V^N(f) = 0$. This estimate extends, modulo some extra technicalities, to Lipschitz functions and then $\|\|\nabla f\|\|_\infty$ is replaced by the Lipschitz norm

$$\|f\|_{\mathcal{L}} := \sup_{x \neq y} \frac{|f(x) - f(y)|}{\|x - y\|},$$

where $x, y$ belong to $\mathbb{R}^{mN^2}$ and $\|x\|$ denotes the Euclidean norm of $x$. Then, for all $\varepsilon > 0$, the following estimate holds:

$$(9) \qquad \mu_V^N(|f - \mu_V^N(f)| > \varepsilon) \leq 2e^{-Nc\varepsilon^2/(2\|f\|_{\mathcal{L}}^2)}.$$

We set

$$f_P(\mathbf{X}) := \hat{\underline{\delta}}^N(P) - m_{P,M}^N$$
$$= \mathrm{Tr}(P(\mathbf{X})) - c_{P,M}^N,$$

with $c_{P,M}^N = \frac{1}{\mu_V^N(\Lambda_M^N)} \int \mathbb{1}_{\Lambda_M^N} \mathrm{Tr}(P(\mathbf{A})) \, d\mu_V^N(\mathbf{A})$. Observing that

$$\partial_{(A_i)_{kl}} \mathrm{Tr}(P(\mathbf{A})) = (D_i P(\mathbf{A}))_{lk},$$

we find that on the closed set $\Lambda_M^N$, $f_P$ is a Lipschitz (actually an infinitely differentiable) function of the entries of $\mathbf{A} \in \mathcal{H}_N(\mathbb{C})^m$ with constant

$$(\|f_P\|_{\mathcal{L}}^{\Lambda_M^N})^2 := \sup_{\mathbf{A} \in \Lambda_M^N} \|\nabla \mathrm{Tr} P(\mathbf{A})\|^2$$

$$= \sup_{\mathbf{A} \in \Lambda_M^N} \sum_{k=1}^m \mathrm{Tr}(D_k P(D_k P)^*) \leq N(\|P\|_{\mathcal{L}}^M)^2,$$

where we simply used that the set of $N \times N$ matrices is a $C^*$-algebra. As a consequence, we also find that, for $\mathbf{B} \in \Lambda_M^N$,

$$|f_P(\mathbf{B})| = \left| \mathrm{Tr}(P(\mathbf{B})) - \frac{1}{\mu_V^N(\Lambda_M^N)} \int \mathbb{1}_{\Lambda_M^N} \mathrm{Tr}(P(\mathbf{A})) \, d\mu_V^N(\mathbf{A}) \right|$$

$$(10) \qquad \leq \frac{\|\mathrm{Tr} P\|_{\mathcal{L}}^{\Lambda_M^N}}{\mu_V^N(\Lambda_M^N)} \int \mathbb{1}_{\Lambda_M^N} \left( \sum_{i=1}^m \mathrm{Tr}(B_i - A_i)^2 \right)^{1/2} d\mu_V^N(\mathbf{A})$$

$$\leq 2\sqrt{m} M N \|P\|_{\mathcal{L}}^M$$

and so on $\Lambda_M^N$ we can extend $f_P$ to $P \in \mathbb{C}\langle X_1, \ldots, X_m \rangle_{\mathcal{L}}^M$.

We can also extend $f_P$ to the whole space $\mathcal{H}_N(\mathbb{C})^m$ with the same Lipschitz constant by putting

$$\bar{f}_P(\mathbf{A}) = \sup_{\mathbf{B} \in \Lambda_M^N} \left\{ f_P(\mathbf{B}) - \sqrt{N} \|P\|_{\mathcal{L}}^M \left( \sum_{i=1}^m \text{Tr}(A_i - B_i)^2 \right)^{1/2} \right\}.$$

Then applying (9), with

$$\varepsilon_{P,M}^N = |\mu_V^N(\mathbb{1}_{(\Lambda_M^N)^c} \bar{f}_P)| + |1 - \mu_V^N(\Lambda_M^N)^{-1}| |\mu_V^N(\mathbb{1}_{\Lambda_M^N} f_P)|,$$

we obtain

$$\begin{aligned}
\mu_V^N(&\{|\hat{\underline{\delta}}^N(P) - m_{P,M}^N| \geq \varepsilon + \varepsilon_{P,M}^N\} \cap \Lambda_M^N) \\
&= \mu_V^N\left( \left\{ \left| \bar{f}_P - \frac{1}{\mu_V^N(\Lambda_M^N)} \mu_V^N(\mathbb{1}_{\Lambda_M^N} \bar{f}_P) \right| \geq \varepsilon + \varepsilon_{P,M}^N \right\} \cap \Lambda_M^N \right) \\
&\leq \mu_V^N(|\bar{f}_P - \mu_V^N(\bar{f}_P)| \geq \varepsilon) \\
&\leq 2e^{-Nc\varepsilon^2/(2(\|f_P\|_{\mathcal{L}}^{\Lambda_M^N})^2)} = 2e^{-c\varepsilon^2/(2(\|P\|_{\mathcal{L}}^M)^2)}.
\end{aligned}$$

We now use the exponential decay of the largest eigenvalue to control $\varepsilon_{P,M}^N$. By (10) and the definition of $\bar{f}_P$, note that

$$\bar{f}_P(\mathbf{A}) \leq 2\sqrt{m} MN \|P\|_{\mathcal{L}}^M + \sqrt{N} \|P\|_{\mathcal{L}}^M \left( \left( \sum_{i=1}^m \text{Tr}(A_i^2) \right)^{1/2} + \sqrt{mN} M \right).$$

Consequently, if $M, N$ are large enough so that $\mu_V^N((\Lambda_M^N)^c) \leq e^{-\alpha NM} \leq \frac{1}{2}$, by Property 2.2,

$$\begin{aligned}
\varepsilon_{P,M}^N &\leq \mu_V^N(\mathbb{1}_{(\Lambda_M^N)^c}(MNm\|P\|_{\mathcal{L}}^M(3 + \lambda_{\max}(\mathbf{A})))) + 2MmN\|P\|_{\mathcal{L}}^M e^{-\alpha NM} \\
&\leq MNm\|P\|_{\mathcal{L}}^M \left( 5e^{-\alpha NM} + \int_0^\infty \mu_V^N(\{\lambda_{\max}(\mathbf{A}) \geq y \vee M\}) \, dy \right) \\
&\leq 6m\left( M^2 + \frac{1}{\alpha N} \right) N\|P\|_{\mathcal{L}}^M e^{-\alpha NM}.
\end{aligned}$$

When $P$ is a monomial of degree $d$,

$$\|P\|_{\mathcal{L}}^M \leq dM^{d-1}.$$

Thus, we only need to control $m_{P,M}^N$;

$$\begin{aligned}
|m_{P,M}^N| &\leq \left| \left( \frac{1}{\mu_V^N(\Lambda_M^N)} - 1 \right) \mu_V^N(\mathbb{1}_{\Lambda_M^N} \text{Tr}(P)) \right| + |\mu_V^N(\mathbb{1}_{(\Lambda_M^N)^c} \text{Tr}(P))| \\
&\leq 2Ne^{-\alpha MN} M^d + N\mu_V^N(\mathbb{1}_{(\Lambda_M^N)^c} \lambda_{\max}(\mathbf{A})^d)
\end{aligned}$$

$$= 2Ne^{-\alpha MN}M^d + dN\int_0^\infty y^{d-1}\mu_V^N(\{\lambda_{\max}(\mathbf{A}) \geq y \vee M\})\, dy$$

$$\leq 2Ne^{-\alpha MN}M^d + dN\int_0^\infty y^{d-1}e^{-\alpha Ny\vee M}\, dy$$

$$\leq (2+1)Ne^{-\alpha MN}M^d + dNe^{-\alpha NM}\sum_{k=1}^d \frac{d-1}{\alpha N}\cdots\frac{d-k}{\alpha N}$$

$$\leq N(3M^d + d^2)e^{-\alpha NM},$$

where we assumed that $d < \alpha N$. $\quad\square$

For later purposes, we have to find a control on the variance of $\hat{\mu}^N$. Recall that $\hat{\underline{\varrho}}^N(P) = N(\hat{\mu}^N(P) - \overline{\mu}^N(P))$.

LEMMA 2.4. *For any* $\varepsilon, \eta, c > 0$, *there exists* $B, C, M_0 > 0$ *such that, for all* $\mathbf{t} \in B_{\eta,c}$, *all* $M \geq M_0$, *for all* $N \in \mathbb{N}$, *and all monomial* $P$ *of degree less than* $\varepsilon N^{(2/3)}$,

$$\mu_V^N((\hat{\underline{\varrho}}^N(P))^2) \leq B(\|P\|_{\mathcal{L}}^M)^2 + C^d N^2 e^{-\alpha MN/2}.$$

PROOF. If $P$ is a monomial of degree $d$, we write

$$\mu_V^N((\hat{\underline{\varrho}}^N(P))^2) \leq \mu_V^N(\mathbb{1}_{\Lambda_M^N}(\hat{\underline{\varrho}}^N(P))^2) + \mu_V^N(\mathbb{1}_{(\Lambda_M^N)^c}(\hat{\underline{\varrho}}^N(P))^2) = I_1 + I_2.$$
(11)

For $I_1$, the previous Lemma implies that

$$I_1 = 2\int_0^\infty x\mu_V^N(\{|\mathrm{Tr}(P) - \mu_V^N(\mathrm{Tr}(P))| \geq x\} \cap \Lambda_M^N)\, dx$$

$$\leq (\varepsilon_{P,M}^N + |m_{P,M}^N|)^2 + 4\int_0^\infty xe^{-cx^2/(2(\|P\|_{\mathcal{L}}^M)^2)}\, dx$$

$$\leq Ce^{-\alpha MN} + B(\|P\|_{\mathcal{L}}^M)^2,$$

with a constant $B = \frac{4}{c}$ and a constant $C$ such that $(\varepsilon_{P,M}^N + |m_{P,M}^N|)^2 \leq Ce^{-\alpha MN}$ for all $d \leq \varepsilon N^{2/3}$. For the second term, we take $M \geq M_0$ with $M_0$ as in Lemma 2.2 (exponential tail of the largest eigenvalue) to get

$$I_2 \leq \mu_V^N[(\Lambda_M^N)^c]^{1/2}\mu_V^N((\hat{\underline{\varrho}}^N(P))^4)^{1/2} \leq e^{-\alpha MN/2}\mu_V^N((\hat{\underline{\varrho}}^N(P))^4)^{1/2}.$$

By the Cauchy–Schwarz inequality, we obtain the control

$$\mu_V^N[\hat{\underline{\varrho}}^N(P)^4] \leq 2^4\mu_V^N((N\hat{\mu}^N(P))^4) \leq 2^4 N^4\mu_V^N((\hat{\mu}^N(PP^*))^2).$$

Now, by the noncommutative Hölder's inequality (see, e.g., [21]),

$$[\hat{\mu}^N(PP^*)]^2 \leq \max_{1\leq i\leq m}\hat{\mu}^N(X_i^{4d})$$

so that we obtain the bound

$$\mu_V^N[\hat{\underline{\varrho}}^N(P)^4] \le 2^4 N^4 \max_{1 \le i \le m} \overline{\mu}^N(X_i^{4d}).$$

By (5) and (6), we obtained a uniform bound $x(= 4\sqrt{Ac^{-1}})$ on $\bar{\mu}^N(X_i)$ so that we have proved using (4) that

$$\overline{\mu}^N(X_i^{4d}) \le 2^{4d}(\mu_c^N(\hat{\mu}^N(X_i^{4d})) + x^{4d}).$$

We can now use the control on the moments as obtained, for instance, by Soshnikov (Theorem 2, page 17 in [22]) to see that there exists $C(\varepsilon)$, $C(\varepsilon) < \infty$ for $\varepsilon > 0$, so that

$$\mu_c^N(\hat{\mu}^N(X_i^{4d})) \le C(\varepsilon)^{4d},$$

provided $d \le \varepsilon N^{2/3}$. As a consequence, we get that

$$\overline{\mu}^N(X_i^{4d}) \le C(\varepsilon)^{4d} \tag{12}$$

for all $d \le \varepsilon N^{2/3}$ and all integer number $N$. Here $C(\varepsilon)$ denotes a finite constant depending only on $\varepsilon$, $\eta$ and $c$ which may have changed from line to line. Hence, we conclude that

$$I_2 \le 4N^2 e^{-\alpha MN/2} C(\varepsilon)^{2d}.$$

Plugging back this estimate into (11), we have proved that for $N$ and $M$ sufficiently large, all monomials $P$ of degree $d \le \varepsilon N^{2/3}$, all $\mathbf{t} \in B_{\eta,c}$,

$$\mu_V^N((\hat{\delta}^N(P))^2) \le B(\|P\|_{\mathcal{L}}^M)^2 + C^{2d}N^2 e^{-\alpha MN/2}$$

with a finite constant $C$ depending only on $\varepsilon$, $c$ and $\eta$. $\square$

## 3. Bound on the distance between $\boldsymbol{\mu}$ and $\overline{\boldsymbol{\mu}}^{\boldsymbol{N}}$.

We here bound, for all monomial $P$,

$$\overline{\delta}^N(P) = N(\overline{\mu}^N(P) - \mu)(P).$$

PROPOSITION 3.1. *For all $c, \varepsilon > 0$, there exists $\eta > 0, C < +\infty$, such that for all integer number $N$, all $\mathbf{t} \in B_{\eta,c}$, and all monomial functions $P$ of degree less than $\varepsilon N^{2/3}$,*

$$|\overline{\delta}^N(P)| \le \frac{C^{deg(P)}}{N}.$$

*In particular, $|(\hat{\underline{\delta}}^N - \hat{\delta}^N)(P)| \le \frac{C^{deg(P)}}{N}$ almost surely.*

Proof. The starting point is the finite dimensional Schwinger–Dyson equation that one gets readily by integration by parts (see [14], proof of Theorem 3.4)

$$(13) \qquad \mu_V^N(\hat{\mu}^N[(X_i + D_iV)P]) = \mu_V^N(\hat{\mu}^N \otimes \hat{\mu}^N(\partial_i P)).$$

Therefore, since $\mu$ satisfies the Schwinger–Dyson equation (2)

$$(14) \qquad \mu[(X_i + D_iV)P] = \mu \otimes \mu(\partial_i P),$$

by taking the difference, we get that for all polynomial $P$,

$$(15) \qquad \overline{\delta}^N(X_i P) = -\overline{\delta}^N(D_iVP) + \overline{\delta}^N \otimes \overline{\mu}^N(\partial_i P) + \mu \otimes \overline{\delta}^N(\partial_i P) + r(N, P),$$

with

$$r(N, P) := N^{-1}\mu_V^N(\hat{\underline{\delta}}^N \otimes \hat{\underline{\delta}}^N(\partial_i P)).$$

If we take $P$ a monomial of degree $d \le \varepsilon N^{2/3}$ and assume $M \ge M_0$, then we see, by using Lemma 2.4,

$$|r(N, P)| \le \frac{1}{N} \sum_{P = P_1 X_i P_2} \mu_V^N(|\hat{\underline{\delta}}^N(P_1)|^2)^{1/2} \mu_V^N(|\hat{\underline{\delta}}^N(P_2)|^2)^{1/2}$$

$$\le \frac{C}{N} \sum_{l=0}^{d-1} (Bl^2 M^{2(l-1)} + C^l N^2 e^{-\alpha MN/2})^{1/2}$$

$$\times (B(d-l-1)^2 M^{2(d-l-1)} + C^{(d-l-1)} N^2 e^{-\alpha MN/2})^{1/2}$$

$$\le \frac{C}{N} d(B(d-1)^2 M^{2(d-2)} + C^{(d-1)} N^2 e^{-\alpha MN/2}) := r(N, d, M).$$

We set

$$\Delta_d^N = \max_{P \text{ monomial of degree } d} |\overline{\delta}^N(P)|.$$

Observe that by (12), for any monomial of degree $d$ less than $\varepsilon N^{2/3}$, $|\overline{\mu}^N(P)| \le C(\varepsilon)^d$, $|\mu(P)| \le C_0^d \le C(\varepsilon)^d$. It allows us to obtain the rough bound $\Delta_d^N \le 2NC(\varepsilon)^d$ if $d < \varepsilon N^{2/3}$. By (15), writing $D_iV = \sum t_j D_i q_j$, we get that, for $d < \varepsilon N^{2/3}$,

$$\Delta_{d+1}^N \le \max_{1 \le i \le m} \sum_{j=1}^n |t_j| \Delta_{d+\deg(D_i q_j)}^N + 2 \sum_{l=0}^{d-1} C(\varepsilon)^{d-l-1} \Delta_l^N + r(N, d, M).$$

We next define, for $\kappa \le 1$,

$$\Delta^N(\kappa, \varepsilon) := \sum_{k=1}^{\varepsilon N^{2/3}} \kappa^k \Delta_k^N.$$

We obtain, if $D$ is the maximal degree of $V$,

$$
(16) \quad
\begin{aligned}
\Delta^N(\kappa, \varepsilon) \leq &[C'\kappa^{-D}|\mathbf{t}| + 2(1 - C(\varepsilon)\kappa)^{-1}\kappa^2]\Delta^N(\kappa, \varepsilon) \\
&+ C|\mathbf{t}| \sum_{k=\varepsilon N^{2/3}+1}^{\varepsilon N^{2/3}+D} \kappa^{k-D}\Delta_k^N + \sum_{k=1}^{\varepsilon N^{2/3}} \kappa^{k+1}r(N, k, M),
\end{aligned}
$$

where we choose $\kappa$ small enough so that $C(\varepsilon)\kappa < 1$. Moreover, since $D$ is finite, using the bound on $\Delta_k^N$, we get

$$
\sum_{k=\varepsilon N^{2/3}+1}^{\varepsilon N^{2/3}+D} \kappa^{k-D}\Delta_k^N \leq 2DN(\kappa C(\varepsilon))^{\varepsilon N^{2/3}}\kappa^{-D}.
$$

Since $\kappa C(\varepsilon) < 1$, as $N$ goes to infinity, this term is negligible with respect to $N^{-1}$ for all $\varepsilon > 0$. The following estimate holds:

$$
\begin{aligned}
&\sum_{k=1}^{\varepsilon N^{2/3}} \kappa^k r(N, k, M) \\
&\leq \frac{C}{N} \sum_{k=1}^{\varepsilon N^{2/3}} k\kappa^k (B(k-1)^2 M^{2(k-2)} + C^{(k-1)}N^2 e^{-\alpha NM/2}) \leq \frac{C''}{N}
\end{aligned}
$$

if $\kappa$ is small enough so that $M^2\kappa < 1$ and $C\kappa < 1$. We observed here that $N^2 e^{-\alpha NM/2}$ is uniformly bounded independently of $N \in \mathbb{N}$. Now, if $|\mathbf{t}|$ is small, we can choose $\kappa$ so that

$$
\zeta := 1 - [C'\kappa^{-D}|\mathbf{t}| + 2(1 - C(\varepsilon)\kappa)^{-1}\kappa^2] > 0.
$$

Plugging these controls into (16) shows that for all $\varepsilon > 0$, and for $\kappa > 0$ small enough, there exists a finite constant $C(\kappa, \varepsilon)$ so that

$$
\Delta^N(\kappa, \varepsilon) \leq C(\kappa, \varepsilon)N^{-1}
$$

and so for all monomial $P$ of degree $d \leq \varepsilon N^{2/3}$,

$$
|\overline{\delta}^N(P)| \leq C(\kappa, \varepsilon)\kappa^{-d}N^{-1}. \qquad \square
$$

To get the precise evaluation of $N\overline{\delta}^N(P)$, we shall first obtain a central limit theorem under $\mu_V^N$ which in turn will allow us to estimate

$$
\lim_{N \to \infty} Nr(N, P).
$$

**4. Central limit theorem.** We shall here prove that

$$\hat{\delta}^N(P) = N(\hat{\mu}^N - \mu)(P)$$

satisfies a central limit theorem for all polynomial $P$. By Proposition 3.1, it is equivalent to prove a central limit theorem for $\underline{\hat{\delta}}^N(P)$, $P \in \mathbb{C}\langle X_1, \ldots, X_m \rangle$. We start by giving a weak form of a central limit theorem for Stieljes-like functions. We then extend the result to polynomial functions in the image of some differential operator. We finally generalize our result to any polynomial functions.

For the rest of the paper, we will always assume the following hypothesis **(H)**.

**(H)**: Let $c$ be a positive real number. The parameter **t** is in $B_{\eta,c}$ with $\eta$ sufficiently small such that we have the convergence to the solution of (2) as well as the control given by Lemma 2.1 (Compact support) and Proposition 3.1.

Note that **(H)** implies also that the control of Lemma 2.1 (Compact support) is uniform, and that we can apply Lemma 2.2 (Exponential tail of the largest eigenvalue) and Lemma 2.3 (Concentration inequality) with uniform constants.

4.1. *Central limit theorem for Stieljes test functions.* One of the issues that one needs to address when working with polynomials is that they are not uniformly bounded. For that reason, we will prefer to work in this section with the complex vector space $\mathcal{C}_{st}^m(\mathbb{C})$ generated by the Stieljes functionals

$$(17) \qquad ST^m(\mathbb{C}) = \left\{ \overrightarrow{\prod_{1 \leq i \leq p}} \left( z_i - \sum_{k=1}^m \alpha_i^k \mathbf{X}_k \right)^{-1} ; \ z_i \in \mathbb{C}\backslash\mathbb{R}, \alpha_i^k \in \mathbb{R}, p \in \mathbb{N} \right\},$$

where $\prod^{\rightarrow}$ is the noncommutative product. We can also equip $ST^m(\mathbb{C})$ with an involution

$$\left( \overrightarrow{\prod_{1 \leq k \leq p}} \left( z_k - \sum_{i=1}^m \alpha_i^k \mathbf{X}_i \right)^{-1} \right)^* = \overrightarrow{\prod_{1 \leq k \leq p}} \left( \overline{z_{p-k}} - \sum_{i=1}^m \alpha_i^{p-k} \mathbf{X}_i \right)^{-1}.$$

We denote $\mathcal{C}_{st}^m(\mathbb{C})_{sa}$ the set of self-adjoint elements of $\mathcal{C}_{st}^m(\mathbb{C})$. The derivation is defined by the Leibniz rule and

$$\partial_i \left( z - \sum_{i=1}^m \alpha_i \mathbf{X}_i \right)^{-1} = \alpha_i \left( z - \sum_{i=1}^m \alpha_i \mathbf{X}_i \right)^{-1} \otimes \left( z - \sum_{i=1}^m \alpha_i \mathbf{X}_i \right)^{-1}.$$

We recall notation; first $\sharp$ is the operator

$$(P \otimes Q)\sharp h = PhQ \quad \text{and} \quad (P \otimes Q \otimes R)\sharp(g, h) = PgQhR$$

so that, for a monomial $q$,

$$\partial_i \circ \partial_j q \#(h_i, h_j) = \sum_{q=q_0 X_i q_1 X_j q_2} q_0 h_i q_1 h_j q_2 + \sum_{q=q_0 X_j q_1 X_i q_2} q_0 h_j q_1 h_i q_2.$$

LEMMA 4.1. *Assume* **(H)** *and let* $h_1, \ldots, h_m$ *be in* $\mathcal{C}_{st}^m(\mathbb{C})_{sa}$. *Then the random variable*

$$Y_N(h_1, \ldots, h_m) = N \sum_{k=1}^m \{\hat{\mu}^N \otimes \hat{\mu}^N(\partial_k h_k) - \hat{\mu}^N[(X_k + D_k V)h_k]\}$$

*converges in law toward a real centered Gaussian variable with covariance*

$$C(h_1, \ldots, h_m) = \sum_{k,l=1}^m (\mu \otimes \mu[\partial_k h_l \times \partial_l h_k] + \mu(\partial_l \circ \partial_k V \sharp(h_k, h_l))) + \sum_{k=1}^m \mu(h_k^2).$$

PROOF. Define $W = \frac{1}{2} \sum_i X_i^2 + V$. Notice that $Y_N(h_1, \ldots, h_m)$ is real valued because the $h_k's$ and $W$ are self adjoint. The proof follows from the usual change of variable trick. We take $h_1, \ldots, h_m$ in $\mathcal{C}_{st}^m(\mathbb{C})_{sa}$, $\lambda \in \mathbb{R}$ and perform a change of variable $A_i \to B_i = F(A)_i = A_i + \frac{\lambda}{N} h_i(\mathbf{A})$ in $Z_V^N$. Note that since the $h_i$ are $\mathcal{C}^\infty$ and uniformly bounded, this defines a bijection on $\mathcal{H}_N(\mathbb{C})^m$ for $N$ big enough. We shall compute the Jacobian of this change of variables up to its second order correction. The Jacobian $J$ may be seen as a matrix $(J_{i,j})_{1 \le i,j \le m}$ where the $J_{i,j}$ are in $\mathcal{L}(\mathcal{H}_N(\mathbb{C}))$ the set of endomorphisms of $\mathcal{H}_N(\mathbb{C})$, and we can write $J = I + \frac{\lambda}{N} \overline{J}$ with

$$\overline{J}_{i,j} : \mathcal{H}_N(\mathbb{C}) \longrightarrow \mathcal{H}_N(\mathbb{C}),$$
$$X \longrightarrow \partial_i h_j \# X.$$

Now, for $1 \le i, j \le m$, $X \longrightarrow \partial_i h_j \# X$ is bounded for the operator norm uniformly in $N$ [since $h_j \in \mathcal{C}_{st}(\mathbb{C}), \partial_i h_j \in \mathcal{C}_{st}(\mathbb{C}) \otimes \mathcal{C}_{st}(\mathbb{C})$ is uniformly bounded] so that, for sufficiently large $N$, the operator norm of $\frac{\lambda}{N} \overline{J}$ is less than 1. From this, we deduce

$$|\det J| = \left|\det\left(I + \frac{\lambda \overline{J}}{N}\right)\right|$$
$$= \exp\left(\operatorname{Tr} \log\left(I + \frac{\lambda \overline{J}}{N}\right)\right) = \exp\left(\sum_{k \ge 1} \frac{(-1)^{k+1} \lambda^k}{k N^k} \operatorname{Tr}(\overline{J}^k)\right).$$

Observe that as $\overline{J}$ is a matrix of size $m^2 N^2$ and of uniformly bounded norm, the $k$th term $\frac{(-1)^{k+1}\lambda^k}{N^k} \operatorname{Tr}(\overline{J}^k)$ is bounded by $\frac{m^2 |\lambda|^k}{N^{k-2}}$. Hence, only the two first terms in the expansion will contribute to the order 1 and the sum $s_N$ of the other terms will be of order $\frac{1}{N}$. To compute the two first terms in the expansion, we only have to remark that if $\phi$ is an endomorphism of $\mathcal{H}_N(\mathbb{C})$ is of the form $\phi(X) = \sum_l A_l X B_l$, with $N \times N$ matrices $A_i, B_i$, then $\mathrm{Tr}\phi = \sum_l \mathrm{Tr}A_l \, \mathrm{Tr}B_l$ [this can be checked by decomposing $\phi$ on the canonical basis of $\mathcal{H}_N(\mathbb{C})$]. Now,

$$\overline{J}_{ij}^k : X \longrightarrow \sum_{1 \le i_1, \ldots, i_{k-1} \le m} \partial_i h_{i_2} \sharp (\partial_{i_2} h_{i_3} \sharp (\cdots (\partial_{i_{k-1}} h_j \sharp X) \cdots)).$$

Thus, we get

$$\mathrm{Tr}(\overline{J}) = \sum_i \mathrm{Tr}\overline{J}_{ii} = \sum_{1 \le i \le m} \mathrm{Tr} \otimes \mathrm{Tr}(\partial_i h_i)$$

and

$$\mathrm{Tr}(\overline{J}^2) = \sum_i \mathrm{Tr}(\overline{J}_{ii}^2) = \sum_{1 \le i, j \le m} \mathrm{Tr} \otimes \mathrm{Tr}(\partial_i h_j \partial_j h_i).$$

We now make the change of variable $A_i \to A_i + \frac{\lambda}{N} h(\mathbf{A})$ to find that

$$Z_V^N = \int e^{-N\mathrm{Tr}(V(\mathbf{A}))} \, d\mu^N(\mathbf{A})$$

$$= \int e^{-N\mathrm{Tr}(W(A_i + (\lambda/N)h_i(\mathbf{A})) - W(A_i))} e^{(\lambda/N) \sum_i \mathrm{Tr} \otimes \mathrm{Tr}(\partial_i h_i)}$$

$$\times e^{(-\lambda^2/(2N^2)) \sum_{i,j} \mathrm{Tr} \otimes \mathrm{Tr}(\partial_i h_j \partial_j h_i)} e^{s_N} \, d\mu_V^N(\mathbf{A})$$

with $s_N$ of order $\frac{1}{N}$. The first term can be expanded into

$$W\left(A_i + \frac{h_i(\mathbf{A})}{N}\right) - W(A_i) = \frac{1}{N} \sum_i \partial_i W \sharp h_i + \frac{1}{N^2} \sum_{i,j} \partial_i \circ \partial_j W \# (h_i, h_j) + \frac{R_N}{N^3},$$

where $R_N$ is a polynomial in the $h_i$'s and in the $X_i$'s, of degree less than the degree of $V$ minus two in the later. To sum up, the following equality holds:

$$\int e^{\lambda Y_N(h_1, \ldots, h_m) - (\lambda^2/2) C_N(h_1, \ldots, h_m) + (1/N)(\hat{\mu}^N(R_N) + N s_N)} = 1,$$

with

$$C_N(h_1, \ldots, h_m) := \hat{\mu}^N\left(\sum_{i,j} \partial_i \circ \partial_j W \# (h_i, h_j)\right) + \hat{\mu}^N \otimes \hat{\mu}^N\left(\sum_{i,j} \partial_i h_j \, \partial_j h_i\right).$$

We can decompose the previous expectation in two terms $E_1$ and $E_2$ with

$$E_1 = \mu_V^N[\mathbb{1}_{\Lambda_M^N} e^{\lambda Y_N(h_1, \ldots, h_m) - (\lambda^2/2) C_N(h_1, \ldots, h_m) + (1/N)(\hat{\mu}^N(R_N) + N s_N)}]$$

and

$$E_2 = \mu_V^N[\mathbb{1}_{(\Lambda_M^N)^c} e^{\lambda Y_N(h_1, \ldots, h_m) - (\lambda^2/2) C_N(h_1, \ldots, h_m) + (1/N)(\hat{\mu}^N(R_N) + N s_N)}].$$

We first consider $E_1$. On $\Lambda_M^N = \{\mathbf{A} : \max_i(\lambda_{\max}^N(A_i)) \le M\}$ the polynomial $R_N$ is uniformly bounded and so $\hat{\mu}^N(R_N) + Ns_N$ is of order one, bounded by a constant $A_N$ which goes uniformly to 0 when $N$ goes to infinity. We next show that we can replace $C_N(h_1, \ldots, h_m)$ by its limit $C(h_1, \ldots, h_m)$ in the exponential in $E_1$. An intermediate step is to replace it by

$$\bar{C}_N(h_1, \ldots, h_m) = \bar{\mu}^N\left(\sum_{i,j} \partial_i \circ \partial_j W \#(h_i, h_j)\right) + \bar{\mu}^N \otimes \bar{\mu}^N\left(\sum_{i,j} \partial_i h_j \, \partial_j h_i\right).$$

In fact, by Lemma 2.3, $\hat{\mu}^N(P)$ converges toward its expectation $\overline{\mu}^N(P)$ under $\mu_V^N(1_{\Lambda_M^N} \cdot)$ except on sets with probability of order $e^{-N^2}$ once evaluated at any products of the $h_i$'s and the $X_i$'s (because the Lipschitz constant of finite products of $h_i$'s and $X_i$'s are bounded on $\Lambda_M^N$ and the error terms $\varepsilon_{P,M}^N$ and $m_{P,M}^N$ can be bounded as we did for polynomials). Hence, we can find a constant $C(M,c) > 0$ such that for $N$ large enough,

$$\mu_V^N(\{|C_N(h_1, \ldots, h_m) - \bar{C}_N(h_1, \ldots, h_m)| > 2\varepsilon\} \cap \Lambda_M^N)$$
$$\le 2e^{-C(M,c)N^2(\varepsilon_N)^2},$$

with $\varepsilon_N = \varepsilon - \varepsilon_{P,M}^N - m_{P,M}^N \sim \varepsilon$. Moreover, $\overline{\mu}^N(P)$ converges to $\mu(P)$ for any polynomial function $P$ (see [14], Theorems 3.1 and 3.4). Since by the Weierstrass theorem the $h_i$'s can be uniformly approximated by polynomials on $\Lambda_M^N$, uniformly in $N$, we also know that $\bar{C}_N(h_1, \ldots, h_m)$ converges to $C(h_1, \ldots, h_m)$. Consequently, we obtain for some positive constant $C(M,c)$, $N$ large enough

$$\mu_V^N(\{|C_N(h_1, \ldots, h_m) - C(h_1, \ldots, h_m)| > \varepsilon\} \cap \Lambda_M^N)$$
$$\le 2e^{-C(M,c)N^2(\varepsilon_N)^2}.$$

Finally, $Y_N(h_1, \ldots, h_m)$ is at most of order $N$ and $C_N(h_1, \ldots, h_m)$ of order one. Hence, the exponential in $E_1$ is at most of order $e^{CN}$ for some finite constant $C$. Therefore, if we let

$$E_1' := \mu_V^N[\mathbb{1}_{\Lambda_M^N} e^{\lambda Y_N(h_1, \ldots, h_m) - (\lambda^2/2)C(h_1, \ldots, h_m)}],$$

we deduce that

$$\left|\log \frac{E_1}{E_1'}\right| \le \left|\log e^{A_N} \frac{\mu_V^N[\mathbb{1}_{\Lambda_M^N} e^{\lambda Y_N(h_1, \ldots, h_m) - (\lambda^2/2)(C_N(h_1, \ldots, h_m) - C(h_1, \ldots, h_m))}]}{\mu_V^N[\mathbb{1}_{\Lambda_M^N} e^{\lambda Y_N(h_1, \ldots, h_m)}]}\right|$$
$$\le |\log(e^{(\lambda^2/2)\varepsilon_N} + 2e^{CN}e^{-C(M,c)N^2\varepsilon_N^2})| + A_N.$$

Letting first $N$ going to infinity and then $\varepsilon$ going to zero yields

$$\lim_{N \to \infty} \frac{E_1}{E_1'} = 1.$$

Note that this estimate is valid for any $M$ large enough so that Lemma 2.3 holds.

Our goal is now to show that, for $M$ sufficiently large, $E_2$ vanishes when $N$ goes to infinity. It would be an easy task if the term in the exponential were bounded, but it may in fact be large due to some derivatives of $V$ appear so that there are polynomials term in the exponential. The idea to pass this difficulty is to make the reverse change of variables. For $N$ bigger than the norm of the $h_i$'s, and with $B_i = A_i + \frac{1}{N} h_i(\mathbf{A})$,

$$E_2 = \mu_V^N[\mathbb{1}_{\{\mathbf{A}:\lambda^N(\mathbf{A})\geq M\}} e^{\lambda Y_N(h_1,\ldots,h_m)-(\lambda^2/2)C_N(h_1,\ldots,h_m)+(1/N)(\hat{\mu}^N(R_N)+Ns_N)}]$$

$$= \mu_V^N(\mathbf{B}:\lambda_{\max}^N(\mathbf{A})\geq M) \leq \mu_V^N(\mathbf{B}:\lambda_{\max}^N(\mathbf{B})\geq M-1).$$

This last quantity goes exponentially fast to 0 for $M$ sufficiently large by Lemma 2.2 (exponential tail of the largest eigenvalue).

Hence, we arrive, for $M$ large enough, at

$$(18) \qquad \lim_{N\to\infty} \int \mathbb{1}_{\Lambda_M^N} e^{\lambda Y_N(h_1,\ldots,h_m)} \, d\mu_V^N = e^{(\lambda^2/2)C(h_1,\ldots,h_m)}.$$

Since $\mu_V^N(\Lambda_M^N)$ goes to one as $N$ goes to infinity, we find that $Y_N(h_1,\ldots,h_m)$ converges in law under $\mu_V^N(\Lambda_M^N)^{-1}\mu_V^N(\cdot \cap \Lambda_M^N)$ toward a centered Gaussian variable with covariance $C(h_1,\ldots,h_m)$, for any $M$ large enough. For the same reason, we conclude that the same convergence holds under $\mu_V^N$. □

4.2. *Central limit theorem for some polynomial functions.* We now extend Lemma 4.1 to polynomial test functions.

LEMMA 4.2. *Assume* (**H**) *and let* $P_1,\ldots,P_m$ *be in* $\mathbb{C}\langle X_1,\ldots,X_m\rangle_{sa}$. *Then, the variable*

$$Y_N(P_1,\ldots,P_m) = N \sum_{k=1}^m [\hat{\mu}^N \otimes \hat{\mu}^N(\partial_k P_k) - \hat{\mu}^N[(X_k+D_kV)P_k]]$$

*converges in law toward a real centered Gaussian variable with covariance*

$$C(P_1,\ldots,P_m) = \sum_{k,l=1}^m (\mu\otimes\mu[\partial_k P_l \times \partial_l P_k] + \mu(\partial_l \circ \partial_k V\sharp(P_k,P_l))) + \sum_{k=1}^m \mu(P_k^2).$$

PROOF. Let $P_1,\ldots,P_m$ be self-adjoint polynomials and $h_1^\varepsilon,\ldots,h_m^\varepsilon$ be Stieljes functionals which approximate $P_1,\ldots,P_m$ such as

$$h_i^\varepsilon(\mathbf{X}) = P_i\left(\frac{X_1}{1+\varepsilon X_1^2},\ldots,\frac{X_m}{1+\varepsilon X_m^2}\right).$$

Since $E[Y_N] = 0$ by (13),

$$Y_N(P_1,\ldots,P_m) = \underline{\hat{\varrho}}^N(K_N(P_1,\ldots,P_m)),$$

with

$$K_N(P_1, \ldots, P_m) = \sum_{k=1}^{m} (\hat{\mu}^N \otimes I(\partial_k P_k) - (X_k + D_k V)P_k)$$

and, similarly, $Y_N(h_1^\varepsilon, \ldots, h_m^\varepsilon) = \hat{\underline{\delta}}^N(K_N(h_1^\varepsilon, \ldots, h_m^\varepsilon))$. It is not hard to see that

$$\|K_N(h_1^\varepsilon, \ldots, h_m^\varepsilon) - K_N(P_1, \ldots, P_m)\|_{\mathcal{L}}^M \leq \varepsilon C(M)$$

for some finite constant $C(M)$ which only depends on $M$. Hence, we deduce by Lemma 2.3 (Concentration inequality) that there exists $m_{P,\varepsilon,M}^N$ and $\varepsilon_{P,\varepsilon,M}^N$ going to zero as $N$ goes to infinity (note here that the control on $m_{P,\varepsilon,M}^N$ and $\varepsilon_{P,\varepsilon,M}^N$ follows exactly the same line as for monomials) such that

$$\mu_V^N(|\hat{\delta}^N(K_N(h_1^\varepsilon, \ldots, h_k^\varepsilon)) - K_N(P_1, \ldots, P_k)) - m_{P,\varepsilon,M}^N| \geq \delta + \varepsilon_{P,\varepsilon,M}^N)$$
$$\leq e^{-\alpha M N} + e^{-\delta^2/(2c\varepsilon^2 C(M)^2)}$$

and so for any bounded continuous function $f: \mathbb{R} \to \mathbb{R}$, if $\nu_{\sigma^2}$ is the centered Gaussian law of covariance $\sigma^2$,

$$\lim_{N \to \infty} \mu_V^N(f(\hat{\underline{\delta}}^N(K_N(P_1, \ldots, P_k)))) = \lim_{\varepsilon \to 0} \lim_{N \to \infty} \mu_V^N(f(\hat{\underline{\delta}}^N(K_N(h_1^\varepsilon, \ldots, h_k^\varepsilon))))$$
$$= \lim_{\varepsilon \to 0} \nu_{C(h_1^\varepsilon, \ldots, h_m^\varepsilon)}(f) = \nu_{C(P_1, \ldots, P_m)}(f),$$

where we used in the second line Lemma 4.2 and in the last line Lemma 2.1 (Compact support) to obtain the convergence of $C(h_1^\varepsilon, \ldots, h_m^\varepsilon)$ to $C(P_1, \ldots, P_m)$.  $\square$

$Y_N$ depends on $N\hat{\mu}^N \otimes \hat{\mu}^N$, in which clearly one of the empirical distribution $\hat{\mu}^N$ shall converge to its deterministic limit. This is the content of the next lemma.

LEMMA 4.3. *Assume* (**H**) *and let* $P_1, \ldots, P_m$ *be self-adjoint polynomial functions. Then, the variable*

$$Z_N(P_1, \ldots, P_m) = \hat{\delta}^N \left( \sum_{k=1}^{m} (X_k + D_k V)P_k - \sum_{k=1}^{m} (\mu \otimes I + I \otimes \mu)(\partial_k P_k) \right)$$

*converges in law toward a centered Gaussian variable with covariance*

$$C(P_1, \ldots, P_m) = \sum_{k,l=1}^{m} (\mu \otimes \mu[\partial_k P_l \times \partial_l P_k] + \mu(\partial_l \circ \partial_k V \sharp(P_k, P_l))) + \sum_{k=1}^{m} \mu(P_k^2).$$

PROOF. The only point is to notice that using (2),

$$Y_N(P_1, \ldots, P_m) = \sum_{k=1}^{m} (\hat{\delta}^N \otimes \mu + \mu \otimes \hat{\delta}^N)(\partial_k P_k) - \hat{\delta}^N((X_k + D_k V)P_k) + r_{N,P}$$

with $r_{N,P} = N^{-1} \sum_{k=1}^{m} \hat{\delta}^N \otimes \hat{\delta}^N(\partial_k P_k)$ of order $N^{-1}$ with probability going to 1 by Lemma 2.3 (Concentration inequality) and Property 3.1. Thus,

$$Y_N(P_1, \ldots, P_m)$$
$$= \hat{\delta}^N \left( \sum_{k=1}^{m} (-(X_k + D_k V)P_k + (I \otimes \mu + \mu \otimes I)(\partial_k P_k)) \right) + r_{N,P}$$
$$= -Z_N(P_1, \ldots, P_m) + O\left(\frac{1}{N}\right).$$

This, with the previous lemma, proves the claim.  □

4.3. *Central limit theorem for all polynomial functions.* In the previous part we have obtained CLT's only for the family of random variables $\hat{\delta}^N(Q)$ with $Q$ in the following subset $\mathcal{F}$ of $\mathbb{C}\langle X_1, \ldots, X_m \rangle$:

$$\mathcal{F} := \left\{ \sum_{k=1}^{m} (X_k + D_k V)P_k - \sum_{k=1}^{m} (\mu \otimes I + I \otimes \mu)(\partial_k P_k), \forall i, P_i \text{ self-adjoint} \right\}.$$

In this section we wish to extend it to $\hat{\delta}^N(Q)$ for any self-adjoint polynomial function $Q$, that is, to prove Theorem 1.3. We have to show a form of density of $\mathcal{F}$ in $\mathbb{C}\langle X_1, \ldots, X_m \rangle$.

The strategy is to see $\mathcal{F}$ as the image of an operator that we will invert. The first operator that comes to mind is

$$\Psi : (P_1, \ldots, P_k) \to \sum_{k=1}^{m} (X_k + D_k V)P_k - \sum_{k=1}^{m} (\mu \otimes I + I \otimes \mu)(\partial_k P_k)$$

as we immediately have $\mathcal{F} = \Psi(\mathbb{C}\langle X_1, \ldots, X_m \rangle_{sa}, \ldots, \mathbb{C}\langle X_1, \ldots, X_m \rangle_{sa})$.

In order to obtain an operator from $\mathbb{C}\langle X_1, \ldots, X_m \rangle_{sa}$ to $\mathbb{C}\langle X_1, \ldots, X_m \rangle_{sa}$, we will prefer to apply $\Psi$ to $P_k = D_k P$ for all $k$ and for a given $P$; as we shall see later, $\Psi(D_1 P, \ldots, D_m P)$ is closely related with the projection on functions of the type $\text{Tr} P$ of the operator on the entries $\Delta - \nabla N \text{Tr}(W).\nabla$ which is symmetric in $L^2(\mu_V^N)$. The resulting operator is a differential operator and, hence, it would be hard to prove that it is continuous on a fixed space of functions. To avoid this issue and make the argument more readable we have first to divide each monomial of $P$ by its degree.

Then, we define a linear map $\Sigma$ on $\mathbb{C}\langle X_1, \ldots, X_m \rangle$ such that, for all monomials $q$ of degree greater or equal to 1,

$$\Sigma q = \frac{q}{\deg q}.$$

Moreover, $\Sigma(q) = 0$ if $\deg q = 0$. For later use, we set $\mathbb{C}_0\langle X_1, \ldots, X_m \rangle$ to be the subset of polynomials $P$ of $\mathbb{C}\langle X_1, \ldots, X_m \rangle_{sa}$ such that $P(0, \ldots, 0) = 0$. We let $\Pi$ be the projection from $\mathbb{C}\langle X_1, \ldots, X_m \rangle_{sa}$ onto $\mathbb{C}_0\langle X_1, \ldots, X_m \rangle$ [i.e., $\Pi(P) = P - P(0, \ldots, 0)$]. We now define some operators on $\mathbb{C}_0\langle X_1, \ldots, X_m \rangle$ that is, from $\mathbb{C}_0\langle X_1, \ldots, X_m \rangle$ into $\mathbb{C}_0\langle X_1, \ldots, X_m \rangle$,

$$\Xi_1 : P \longrightarrow \Pi\left( \sum_{k=1}^m \partial_k \Sigma P \sharp D_k V \right),$$

$$\Xi_2 : P \longrightarrow \Pi\left( \sum_{k=1}^m (\mu \otimes I + I \otimes \mu)(\partial_k D_k \Sigma P) \right).$$

We denote $\Xi_0 = Id - \Xi_2$ and $\Xi = \Xi_0 + \Xi_1$, where $I$ is the identity on $\mathbb{C}_0\langle X_1, \ldots, X_m \rangle$. Note that the images $\Xi_i$'s and $\Xi$ are indeed included in $\mathbb{C}\langle X_1, \ldots, X_m \rangle_{sa}$ since $V$ is assumed self-adjoint. With this notation, Lemma 4.3, once applied to $P_i = D_i \Sigma P$, $1 \leq i \leq m$, reads as follows:

PROPOSITION 4.4. *For all $P$ in $\mathbb{C}_0\langle X_1, \ldots, X_m \rangle$, $\hat{\delta}^N(\Xi P)$ converges in law to a centered Gaussian variable with covariance*

$$\mathcal{C}(P) := C(D_1 \Sigma P, \ldots, D_m \Sigma P).$$

PROOF. We have for all tracial state $\tau$, $\tau(\partial_k P \sharp V) = \tau(D_k P V)$ and if $P$ is in $\mathbb{C}_0\langle X_1, \ldots, X_m \rangle$ [i.e., $P(0, \ldots, 0) = 0$], we have the identity

$$P = \sum_k \partial_k \Sigma P \sharp X_k.$$

Then, as $\hat{\delta}^N$ is tracial and null on constant terms (so that the projection $\Pi$ can be removed in the definition of $\Xi$), for all polynomial $P$,

$$\hat{\delta}^N(\Xi P) = \hat{\delta}^N\left( P + \sum_{k=1}^m \partial_k \Sigma P \sharp D_k V - \sum_{k=1}^m (\mu \otimes I + I \otimes \mu)(\partial_k D_k \Sigma P) \right)$$

$$= \hat{\delta}^N\left( \sum_{k=1}^m (X_k + D_k V) D_k \Sigma P - \sum_{k=1}^m (\mu \otimes I + I \otimes \mu)(\partial_k D_k \Sigma P) \right)$$

$$= Z_N(D_1 \Sigma P, \ldots, D_m \Sigma P).$$

We then use the Lemma 4.3 to conclude. $\quad\square$

To generalize the central limit theorem to all polynomial functions, we need to show that the image of $\Xi$ is dense and to control approximations. If $P$ is a polynomial and $q$ a nonconstant monomial, we will denote $\lambda_q(P)$ the coefficient of $q$ in the decomposition of $P$ in monomials. We can then define a norm $\| \cdot \|_A$ on $\mathbb{C}_0\langle X_1, \ldots, X_m \rangle$ for $A > 1$ by

$$\|P\|_A = \sum_{\deg q \neq 0} |\lambda_q(P)| A^{\deg q}.$$

In the formula above, the sum is taken on all nonconstant monomials. We also define the operator norm given, for $T$ from $\mathbb{C}_0\langle X_1, \ldots, X_m \rangle$ to $\mathbb{C}_0\langle X_1, \ldots, X_m \rangle$, by

$$\|T\|_A = \sup_{\|P\|_A = 1} \|T(P)\|_A.$$

Finally, let $\mathbb{C}_0\langle X_1, \ldots, X_m \rangle_A$ be the completion of $\mathbb{C}_0\langle X_1, \ldots, X_m \rangle$ for $\| \cdot \|_A$. We say that $T$ is continuous on $\mathbb{C}_0\langle X_1, \ldots, X_m \rangle_A$ if $\|T\|_A$ is finite. We shall prove that $\Xi$ is continuous on $\mathbb{C}_0\langle X_1, \ldots, X_m \rangle_A$ with continuous inverse when $\mathbf{t}$ is small.

Lemma 4.5.  *With the previous notation:*

1. *The operator $\Xi_0$ is invertible on $\mathbb{C}_0\langle X_1, \ldots, X_m \rangle$.*
2. *There exists $A_0 > 0$ such that, for all $A > A_0$, the operators $\Xi_2$, $\Xi_0$ and $\Xi_0^{-1}$ are continuous on $\mathbb{C}_0\langle X_1, \ldots, X_m \rangle_A$ and their norms are uniformly bounded for $\mathbf{t}$ in $B_\eta$.*
3. *For all $\varepsilon, A > 0$, there exists $\eta_\varepsilon > 0$ such for $|\mathbf{t}| < \eta_\varepsilon$, $\Xi_1$ is continuous on $\mathbb{C}_0\langle X_1, \ldots, X_m \rangle_A$ and $\|\Xi_1\|_A \leq \varepsilon$.*
4. *For all $A > A_0$, there exists $\eta > 0$ such that for $\mathbf{t} \in B_\eta$, $\Xi$ is continuous, invertible with a continuous inverse on $\mathbb{C}_0\langle X_1, \ldots, X_m \rangle_A$. Besides, the norms of $\Xi$ and $\Xi^{-1}$ are uniformly bounded for $\mathbf{t}$ in $B_\eta$.*
5. *There exists $C > 0$ such that, for all $A > C$, $\mathcal{C}$ is continuous from $\mathbb{C}_0\langle X_1, \cdots, X_m \rangle_A$ into $\mathbb{R}$.*

Proof.  1. We can write

$$\Xi_0 = I - \Xi_2.$$

Observe that since $\Xi_2$ reduces the degree of a polynomial by at least 2,

$$P \to \sum_{n \geq 0} (\Xi_2)^n(P)$$

is well defined on $\mathbb{C}_0\langle X_1, \ldots, X_m \rangle$ as the sum is finite for any polynomial $P$. This clearly gives an inverse for $\Xi_0$.

2. First remark that a linear operator $T$ has a norm less than $C$ with respect to $\| \cdot \|_A$ if and only if for all nonconstant monomial $q$,

$$\|T(q)\|_A \leq C A^{\deg q}.$$

Recall that $\mu$ is uniformly compactly supported [see Lemma 2.1 (Compact support)] and let $C_0 < +\infty$ be such that $|\mu(q)| \leq C_0^{\deg q}$ for all monomial $q$. Take a monomial $q = X_{i_1} \cdots X_{i_p}$, and assume that $A > 2C_0$,

$$\left\| \Pi\left( \sum_k (I \otimes \mu)\partial_k D_k \Sigma q \right) \right\|_A$$

$$\leq p^{-1} \sum_{\substack{k, q = q_1 X_k q_2, \\ q_2 q_1 = r_1 X_k r_2}} \| r_1 \mu(r_2) \|_A$$

$$\leq p^{-1} \sum_{\substack{k, q = q_1 X_k q_2, \\ q_2 q_1 = r_1 X_k r_2}} A^{\deg r_1} C_0^{\deg r_2} = \frac{1}{p} \sum_{n=0}^{p-1} \sum_{l=0}^{p-2} A^l C_0^{p-l-2}$$

$$\leq A^{p-2} \sum_{l=0}^{p-2} \left( \frac{C_0}{A} \right)^{p-2-l} \leq 2A^{-2} \| q \|_A,$$

where in the second line, we observed that once $\deg(q_1)$ is fixed, $q_2 q_1$ is uniquely determined and then $r_1, r_2$ are uniquely determined by the choice of $l$ the degree of $r_1$. Thus, the factor $\frac{1}{p}$ is compensated by the number of possible decompositions of $q$, that is, the choice of $n$ the degree of $q_1$. If $A > 2$, $P \to \Pi(\sum_k (I \otimes \mu)\partial_k D_k \Sigma P)$ is continuous of norm strictly less than $\frac{1}{2}$. And a similar calculus for $\Pi(\sum_k (\mu \otimes I)\partial_k D_k \Sigma)$ shows that $\Xi_2$ is continuous of norm strictly less than 1. It follows immediately that $\Xi_0$ is continuous. Recall now that

$$\Xi_0^{-1} = \sum_{n \geq 0} \Xi_2^n.$$

As $\Xi_2$ is of norm strictly less than 1, $\Xi_0^{-1}$ is immediately continuous.

3. Let $q = X_{i_1} \cdots X_{i_p}$ be a monomial and let $D$ be the degree of $V$ and $B(\leq Dn)$ the sum of the maximum number of monomials in $D_k V$:

$$\| \Xi_1(q) \|_A \leq \frac{1}{p} \sum_{k, q = q_1 X_k q_2} \| q_1 D_k V q_2 \|_A$$

$$\leq \frac{1}{p} \sum_{k, q = q_1 X_k q_2} |\mathbf{t}| B A^{p-1+D-1}$$

$$= |\mathbf{t}| B A^{D-2} \| q \|_A.$$

It is now sufficient to take $\eta_\varepsilon < (BA^{D-2})^{-1}\varepsilon$.

4. We choose $\eta < (BA^{D-2})^{-1} \| \Xi_0^{-1} \|_A^{-1}$ so that when $|\mathbf{t}| \leq \eta$,

$$\| \Xi_1 \|_A \| \Xi_0^{-1} \|_A < 1.$$

By continuity, we can extend $\Xi_0$, $\Xi_1$, $\Xi_2$, $\Xi$ and $\Xi_0^{-1}$ on the space $\mathbb{C}_0\langle X_1, \ldots, X_m \rangle_A$. The operator

$$P \to \sum_{n \geq 0} (-\Xi_0^{-1}\Xi_1)^n \Xi_0^{-1}$$

is well defined and continuous. And this is clearly an inverse of

$$\Xi = \Xi_0 + \Xi_1 = \Xi_0(I + \Xi_0^{-1}\Xi_1).$$

5. We finally prove that $\mathcal{C}$ is continuous from $\mathbb{C}_0\langle X_1, \ldots, X_m \rangle_A$ into $\mathbb{R}$ where we recall that we assumed $A > C_0$. Let us consider the first term

$$\mathcal{C}_1(P) := \sum_{k,l=1}^{m} \mu \otimes \mu(\partial_k D_l \Sigma P \times \partial_l D_k \Sigma P).$$

Then we obtain, as in the second point of this proof,

$$|\mathcal{C}_1(P)| \leq 4 \sum_{k,l=1}^{m} \sum_{q,q'} \frac{|\lambda_q(P)||\lambda_{q'}(P)|}{\deg q \deg q'} \sum_{\substack{q=q_1 X_k q_2, q'=q_1' X_l q_2' \\ q_2 q_1 = r_1 X_l r_2, q_2' q_1' = r_1' X_k r_2'}} C_0^{\deg q + \deg q' - 4}$$

$$\leq 4 \sum_{q,q'} |\lambda_q(P)||\lambda_{q'}(P)| \deg q \deg q' C_0^{\deg q + \deg q' - 4}$$

$$\leq 4 \left( \sup_{\ell \geq 0} \ell C_0^{\ell-2} A^{-\ell} \right)^2 \|P\|_A^2.$$

We next turn to show that

$$\mathcal{C}_2(P) := \sum_{k,l=1}^{m} \mu(\partial_k \circ \partial_l V \sharp (D_k \Sigma P, D_l \Sigma P))$$

is also continuous for $\|\cdot\|_A$. In fact, noting that we may assume $V \in \mathbb{C}_0\langle X_1, \ldots, X_m \rangle$ without changing $\mathcal{C}_2$,

$$|\mathcal{C}_2(P)| \leq \sum_{p,q,q',k,l} |\lambda_p(V)|$$

$$\times \sum_{\substack{q,q',p=p_1 X_k p_2 X_l p_3 \\ q=q_1 X_k q_2, q'=q_1' X_k q_2'}} \frac{|\lambda_q(P)||\lambda_{q'}(P)| C_0^{\deg p + \deg q + \deg q' - 4}}{\deg q \deg q'}$$

$$\leq n|\mathbf{t}|D^2 \sum_{q,q'} |\lambda_q(P)||\lambda_{q'}(P)| C_0^{D + \deg q + \deg q' - 4}$$

$$\leq n|\mathbf{t}|D^2 C_0^{D-4} \|P\|_A^2.$$

The continuity of the last term $\mathcal{C}_3(P) = \sum_{i=1}^{m} \mu((D_j \Sigma P)^2)$ is obtained similarly. $\square$

We can compare the norm $\| \cdot \|_A$ to a more intuitive norm, namely, $\| \cdot \|_{\mathcal{L}}^M$ defined in (7).

We will say that a semi-norm $\mathcal{N}$ is weaker than a semi-norm $\mathcal{N}'$ if and only if there exists $C < +\infty$ such that, for all $P$ in $\mathbb{C}_0\langle X_1, \ldots, X_m \rangle$,

$$\mathcal{N}(P) \leq C\mathcal{N}'(P).$$

LEMMA 4.6. *For $A > M$, the semi-norm $\| \cdot \|_{\mathcal{L}}^M$ restricted to the space $\mathbb{C}_0\langle X_1, \ldots, X_m \rangle$ is weaker than the norm $\| \cdot \|_A$.*

PROOF. For all $P$ in $\mathbb{C}_0\langle X_1, \ldots, X_m \rangle$, the following inequalities hold:

$$\|P\|_{\mathcal{L}}^M \leq \sum_q |\lambda_q(P)| \|q\|_{\mathcal{L}}^M \leq \sum_q |\lambda_q(P)| \deg q M^{\deg q} \leq \left( \sup_l l \left( \frac{M}{A} \right)^l \right) \|P\|_A. \quad \square$$

To take into account the previous results, we define a new hypothesis ($\mathbf{H}'$) stronger than ($\mathbf{H}$).

($\mathbf{H}'$): ($\mathbf{H}$) is satisfied, $A - 1 > \max(A_0, M_0, C)$ for the $M_0$ which appear in Lemma 2.2 (Exponential tail of the largest eigenvalue) and the $C$ which appear in Proposition 3.1. Besides, $|\mathbf{t}| \leq \eta$ with $\eta$ as in the fourth point of Lemma 4.5 in order that $\Xi$ and $\Xi^{-1}$ are continuous on $\mathbb{C}_0\langle X_1, \ldots, X_m \rangle_A$ and $\mathbb{C}_0\langle X_1, \ldots, X_m \rangle_{A-1}$, and that $\mathcal{C}$ is also continuous for these norms.

The two main additional consequences of this hypothesis are the continuity of $\Xi$ for $\| \cdot \|_A$. The strange condition about the continuity of $\Xi$ on $\mathbb{C}_0\langle X_1, \ldots, X_m \rangle_{A-1}$ is here for a technical reason which will appear only in the last section on the interpretation of the first order correction to the free energy.

While ($\mathbf{H}'$) is full of conditions, the only important hypothesis is the $c$-convexity of $V$. Given such a $V$, we can always find constants $A$ and $\eta$ which satisfy the hypothesis. The only restriction will be then that $\mathbf{t}$ is sufficiently small.

We can now prove the general central limit theorem which is up to the identification of the covariance equivalent to Theorem 1.3.

THEOREM 4.7. *Assume ($\mathbf{H}'$). For all $P$ in $\mathbb{C}\langle X_1, \ldots, X_m \rangle_{sa}$, $\hat{\delta}^N(P)$ converges in law to a centered Gaussian variable $\gamma_P$ with covariance*

$$\sigma^2(P) := \mathcal{C}(\Xi^{-1}\Pi(P)) = C(D_1 \Sigma \Xi^{-1}\Pi(P), \ldots, D_m \Sigma \Xi^{-1}\Pi(P)).$$

*If $P \in \mathbb{C}\langle X_1, \ldots, X_m \rangle$, $\hat{\delta}^N(P)$ converges to the complex centered Gaussian variable $\gamma_{(P+P^*)/2} + i\gamma_{(P-P^*)/2i}$ [the covariance of $\gamma_{(P+P^*)/2}$ and $\gamma_{(P-P^*)/2i}$ being given by $\sigma^2((P+P^*)/2, (P-P^*)/2i)$, where $\sigma^2(\cdot, \cdot)$ is the bilinear form associated to the quadratic form $\sigma^2$].*

PROOF. As $\hat{\delta}^N(P)$ does not depend on constant terms, we can directly take $P = \Pi(P)$ in $\mathbb{C}_0\langle X_1, \ldots, X_m \rangle$. Now, by part 4 of Lemma 4.5, we can find an element $Q$ of $\mathbb{C}_0\langle X_1, \ldots, X_m \rangle_A$ such that $\Xi Q = P$. But the space $\mathbb{C}_0\langle X_1, \ldots, X_m \rangle$ is dense in $\mathbb{C}_0\langle X_1, \ldots, X_m \rangle_A$ by construction. Thus, there exists a sequence $Q_n$ in $\mathbb{C}_0\langle X_1, \ldots, X_m \rangle$ such that

$$\lim_{n \to \infty} \|Q - Q_n\|_A = 0.$$

Let us define $R_n = P - \Xi Q_n$ in $\mathbb{C}_0\langle X_1, \ldots, X_m \rangle$.

Now according to Property 4.4 for all $n$, $\hat{\delta}^N(\Xi Q_n)$ converges in law to a Gaussian variable $\gamma_n$ of variance $\mathcal{C}(Q_n)$ with

$$\mathcal{C}(Q_n) = C(D_1 \Sigma Q_n, \ldots, D_m \Sigma Q_n).$$

As $\mathcal{C}$ is continuous by part 4 of Lemma 4.5, it can be extended to the space $\mathbb{C}_0\langle X_1, \ldots, X_m \rangle_A$ and $\sigma^2(P) = \mathcal{C}(\Xi^{-1}P) = \mathcal{C}(Q) = \lim_n \mathcal{C}(Q_n)$ is well defined. Hence, $\gamma_n$ converges weakly toward $\gamma_\infty$, the centered Gaussian law with covariance $\mathcal{C}(Q)$, when $n$ goes to $+\infty$. The last step is to prove the convergence in law of $\hat{\delta}^N(P)$ to $\gamma_\infty$. We will use the Dudley distance. For $f : \mathbb{R} \to \mathbb{R}$, we define $|f|_{\mathcal{L}} = \|f\|_{\mathcal{L}} + \|f\|_\infty$. The Dudley distance between two measures on $\mathbb{R}$ is

$$\mathcal{D}(\mu, \nu) = \sup_{|f|_{\mathcal{L}} \leq 1} |\mu(f) - \nu(f)|.$$

The topology induced by the Dudley metric is the topology of the convergence in law. Below, as a parameter of $\mathcal{D}$, we denote in short $\hat{\delta}^N(P)$ for the law of $\hat{\delta}^N(P)$. We make the following decomposition:

(19)
$$\begin{aligned} \mathcal{D}(\hat{\delta}^N(P), \gamma_\infty) \leq{} &\mathcal{D}(\hat{\delta}^N(P), \hat{\delta}^N(\Xi Q_n)) \\ &+ \mathcal{D}(\hat{\delta}^N(\Xi Q_n), \gamma_n) + \mathcal{D}(\gamma_n, \gamma_\infty). \end{aligned}$$

By the above remarks, $\mathcal{D}(\hat{\delta}^N(\Xi Q_n), \gamma_n)$ goes to 0 when $N$ goes to $+\infty$ and $\mathcal{D}(\gamma_n, \gamma_\infty)$ goes to 0 when $n$ goes to $+\infty$. We now use the bound on the Dudley distance:

$$\mathcal{D}(\hat{\delta}^N(P), \hat{\delta}^N(\Xi Q_n)) \leq E[|\hat{\delta}^N(P) - \hat{\delta}^N(\Xi Q_n)| \wedge 1] = E[|\hat{\delta}^N(R_n)| \wedge 1].$$

We control the last term by Lemmas 2.3 (Concentration inequality) and 2.2 (Exponential tail of the largest eigenvalue) so that, for $M \geq M_0$,

$$E[|\hat{\delta}^N(R_n)| \wedge 1] \leq e^{-\alpha NM} + 2\sqrt{\frac{2\pi}{c}} \|R_n\|_{\mathcal{L}}^M + \varepsilon_{R_n, M}^N + |m_{R_n, M}^N|.$$

But we deduce from Lemma 4.6 that since we chose $M < A$, there exists a finite constant $C$ such that

$$\|R_n\|_{\mathcal{L}}^M \leq C\|R_n\|_A = C\|\Xi(Q - Q_n)\|_A \leq C\|\Xi\|_A\|Q - Q_n\|_A$$

and so $\|R_n\|_{\mathcal{L}}^M$ goes to zero as $n$ goes to infinity. And since $\|R_n\|_{\mathcal{L}}^M$ is finite, $\varepsilon_{R_n,M}^N$ goes to zero. Similarly, using the bound of Lemma 2.3 on $m_{P,M}^N$ for $P$ monomial, we find that

$$|m_{R_n,M}^N| \leq N \sum_q |\lambda_q(R_n)| \deg(q)(3M^{\deg(q)} + \deg(q)^2)e^{-\alpha MN}$$

$$\leq N \sup_{\ell \geq 0}(\ell(3M^\ell + \ell^2)A^{-\ell})\|R_n\|_A e^{-\alpha MN}$$

goes to zero as $N$ goes to infinity. Thus, $E[|\hat{\delta}^N(R_n)| \wedge 1]$ goes to zero as $n$ and $N$ go to infinity. Putting things together, we obtain if we let first $N$ going to $+\infty$ and then $n$, the desired convergence $\lim_N \mathcal{D}(\hat{\delta}^N(P), \gamma_\infty) = 0$. $\square$

Note that the convergence in law in Theorem 4.7 can be generalized to a convergence *in moments*;

COROLLARY 4.8. *Assume* ($\mathbf{H}'$). *Let $P$ be a self-adjoint polynomial, then $\hat{\delta}^N(P)$ converges in moments to a real centered Gaussian variable with variance $\sigma^2(P)$, that is, for all $k$ in $\mathbb{N}$,*

$$\lim_{N \to \infty} \int (\hat{\delta}^N P)^k \, d\mu_V^N = \frac{1}{\sqrt{2\pi\sigma^2(P)}} \int x^k e^{-x^2/(2\sigma^2(P))} \, dx.$$

PROOF. Indeed, once again we decompose $\int (\hat{\delta}^N P)^k \, d\mu_V^N$ into $E_1^N + E_2^N$ with

$$E_1^N = \int \mathbb{1}_{\Lambda_M^N}(\hat{\delta}^N P)^k \, d\mu_V^N \qquad E_2^N = \int \mathbb{1}_{(\Lambda_M^N)^c}(\hat{\delta}^N P)^k \, d\mu_V^N,$$

with $M \geq M_0$. For $E_1$, we notice that the law of $\hat{\delta}^N P$ has a sub-Gaussian tail according to Lemma 2.3 (Concentration inequality). Therefore, we can replace $x^k$ by a bounded continuous function, producing an error independent of $N$. Applying Theorem 4.7 then shows that

$$\lim_{N \to \infty} \int \mathbb{1}_{\Lambda_M^N}(\hat{\delta}^N P)^k \, d\mu_V^N = \frac{1}{\sqrt{2\pi\sigma^2(P)}} \int x^k e^{-x^2/(2\sigma^2(P))} \, dx.$$

For the second term, we use the trivial bound

$$|E_2^N| \leq N^k \int \mathbb{1}_{(\Lambda_M^N)^c}(|\lambda_{\max}(\mathbf{A})| + |\mu|(P))^k \, d\mu_V^N$$

$$\leq kN^k \int_{\lambda \geq M}(\lambda + |\mu|(P))^{k-1}e^{-\alpha\lambda N} \, d\lambda,$$

which goes to zero as $N$ goes to infinity for all finite $k$. $\square$

Another generalization of Theorem 4.7 is to extend the set of test functions from polynomials to the completion of $\mathbb{C}_0\langle X_1,\ldots,X_m\rangle$ for the Lipschitz semi-norm $\|\cdot\|_{\mathcal{L}}^M$. We shall assume that $M$ is strictly greater than $C$, the constant which bounds uniformly the radius of the support of $\mu$ according to Lemma 2.1 (Compact support), and also greater than $M_0$, the constant which appears in Lemma 2.2 (Exponential tail of the largest eigenvalue) in order to have $\lambda_{\max}(\mathbf{A})$ less than $M$ with high probability. We denote $\mathbb{C}_0\langle X_1,\ldots,X_m\rangle_{\mathcal{L}}^M$ the completion of $\mathbb{C}_0\langle X_1,\ldots,X_m\rangle$ for that norm.

Let us first extend some of the previous quantities to this setting. Recall that, for all $N\in\mathbb{N}$, $\sqrt{N}\|P\|_{\mathcal{L}}^M$ is always bigger than $\|\mathrm{Tr}P\|_{\mathcal{L}}^{\Lambda_M^N}$, so that if $\lambda_{\max}(\mathbf{A})<M$, $\mathrm{Tr}P(\mathbf{A})$ is well defined. This allows us to define, for $P$ in $\mathbb{C}_0\langle X_1,\ldots,X_m\rangle_{\mathcal{L}}^M$, $\hat{\mu}^N(P)=\frac{1}{N}\mathrm{Tr}P(\mathbf{A})$ on $\Lambda_M^N$. We can also extend $\mu$ to this context by the following:

LEMMA 4.9.　*Let $P\in\mathbb{C}_0\langle X_1,\ldots,X_m\rangle$. Then, with $C_0$ as in Lemma 2.1 (Compact support),*

$$|\mu(P)|\leq\sqrt{m}C_0\|P\|_{\mathcal{L}}^{C_0}.$$

PROOF.　Let us consider the following norm on $\mathbb{C}\langle X_1,\ldots,X_m\rangle$:

$$\|P\|_\mu:=\limsup_n(\mu((PP^*)^n))^{1/(2n)}.$$

The completion and separation of $\mathbb{C}\langle X_1,\ldots,X_m\rangle$ for this norm is then a $C^*$-algebra (see, e.g., the Gelfand–Neimark–Segal construction). As $\mu$ is compactly supported, the norm of the $X_i$'s are bounded by $C_0$. Besides, for all $P$,

$$|\mu(P)|\leq\|P\|_\mu.$$

Therefore, we can write

$$|\mu(P)|=|\mu(P(\mathbf{X}))-\mu(P(0))|=\left|\mu\left(\int_0^1\sum_{k=1}^m(D_kP)(s\mathbf{X})X_k\,ds\right)\right|$$

$$\leq\int_0^1\left|\mu\left(\sum_{k=1}^mD_kP(s\mathbf{X})X_k\right)\right|ds$$

$$\leq\int_0^1\left(\sum_{k=1}^m\mu(D_kP(s\mathbf{X})D_kP(s\mathbf{X})^*)\right)^{1/2}\left(\sum_{k=1}^m\mu(X_k^2)\right)^{1/2}ds$$

$$\leq C_0\sup_{\substack{\mathcal{A}C^*\text{-algebra}\\x_i=x_i^*\|x_i\|\leq C_0}}\left(\sum_{k=1}^m\|D_kP(x_1,\ldots,x_m)\|_{\mathcal{A}}^2\right)^{1/2}=\sqrt{m}C_0\|P\|_{\mathcal{L}}^{C_0}.\qquad\square$$

Thus, $\mu$ extends to $\mathbb{C}\langle X_1, \ldots, X_m \rangle_{\mathcal{L}}^M$. It is a natural question to study the behavior of

$$\hat{\delta}^N(P) := N(\hat{\mu}^N(P) - \mu(P)) 1_{\Lambda_M^N}$$

for $P$ in $\mathbb{C}_0\langle X_1, \ldots, X_m \rangle_{\mathcal{L}}^M$, the completion of $\mathbb{C}_0\langle X_1, \ldots, X_m \rangle$ for $\| \cdot \|_{\mathcal{L}}^M$.

COROLLARY 4.10. *Assume* (**H**$'$) *and let* $M$ *be bigger than* $C_0$ *and* $M_0$:

1. $\sigma^2$ *is continuous for* $\| \cdot \|_{\mathcal{L}}^M$ *and so extends to* $\mathbb{C}_0\langle X_1, \ldots, X_m \rangle_{\mathcal{L}}^M$.
2. *For all* $P$ *in* $\mathbb{C}_0\langle X_1, \ldots, X_m \rangle_{\mathcal{L}}^M$, $\hat{\delta}^N(P)$ *converges in law to a Gaussian variable with variance* $\sigma^2(P)$.

PROOF. We take a sequence of polynomials $S_n$ which converges to $P$ for the norm $\| \cdot \|_{\mathcal{L}}^M$. Let $R_n = P - S_n$ be the rest. For all $n$, $\hat{\delta}^N(S_n)$ converges to a centered Gaussian variable $\gamma_n$ of variance $\sigma^2(S_n)$.

Let us show that $\sigma^2$ is continuous for $\| \cdot \|_{\mathcal{L}}^M$. Let $P$ be a polynomial, and $M$ sufficiently large,

$$\sigma^2(P) = \lim_N E[\hat{\underline{\delta}}^N(P)^2] = \lim_N E[1_{\Lambda_M^N} \hat{\underline{\delta}}^N(P)^2].$$

The first equality comes from the previous corollary about the convergence in moments, as well as Lemma 3.1, which allows to recenter with respect to the mean rather than the limit, and the second equality comes from Lemma 2.2 (Exponential tail of the largest eigenvalue). Now by Lemma 2.3 (Concentration inequality), as $\|P\|_{\mathcal{L}}^M$ controls the Lipschitz norm of $\frac{1}{N}\text{Tr}(P)$,

$$\lim_{N \to \infty} \mu_V^N[1_{\Lambda_M^N}(\hat{\underline{\delta}}^N(P))^2]$$

$$= 2 \lim_{N \to \infty} \int_0^\infty \varepsilon \mu_V^N(\Lambda_M^N \cap \{|\hat{\underline{\delta}}^N(P)| > \varepsilon\}) \, d\varepsilon$$

$$\leq \int_0^\infty 2\varepsilon e^{-c\varepsilon^2/(2(\|P\|_{\mathcal{L}}^M)^2)} \, d\varepsilon = \frac{4}{c}(\|P\|_{\mathcal{L}}^M)^2,$$

where we used that $m_{M,P}^N$ and $\varepsilon_{M,P}^N$ of Lemma 2.3 go to zero as $N$ goes to infinity since $P$ is a polynomial. Thus, the quadratic form $\sigma^2$ is continuous for $\| \cdot \|_{\mathcal{L}}^M$ and can be extended on $\mathbb{C}_0\langle X_1, \ldots, X_m \rangle_{\mathcal{L}}^M$. This implies that $\sigma^2(S_n)$ converges to $\sigma^2(P)$. The rest of the proof is exactly as that of Theorem 4.7 and we omit it. $\square$

Note that by Lemma 4.5 the norm $\| \cdot \|_A$ is stronger than the norm $\| \cdot \|_{\mathcal{L}}^M$ so that we can use this corollary to extend out the central limit theorem on $\mathbb{C}_0\langle X_1, \ldots, X_m \rangle_A$ and, by continuity of $\sigma^2$, on this space the formula

$$\sigma^2(P) := \mathcal{C}(\Xi^{-1}P) = C(D_1\Sigma\Xi^{-1}P, \ldots, D_m\Sigma\Xi^{-1}P)$$

remains valid.

**5. Identification of the variance.**

5.1. *Exact formula.* We shall provide here a more tractable formula for
the variance $\sigma^2(P)$ of the limiting Gaussian distribution found in Theorem
4.7. Note that for all polynomials $P$, $Q$, $\hat{\delta}^N(P+Q)$ converges to $\gamma_{P+Q}$.
Thus, $\{\gamma_P | P \in \mathbb{C}\langle X_1, \ldots, X_m \rangle_{sa}\} = \{\gamma_P | P \in \mathbb{C}_0\langle X_1, \ldots, X_m \rangle\}$ has a
natural structure of Gaussian space. In this space all elements are centered and the
covariance function is given, for $P, Q \in \mathbb{C}_0\langle X_1, \ldots, X_m \rangle$ by

$$\sigma^2(P, Q) = \mathcal{C}(\Xi^{-1}P, \Xi^{-1}Q) = C(\mathbf{D}\Sigma\Xi^{-1}P, \mathbf{D}\Sigma\Xi^{-1}Q),$$

where $\mathbf{D}$ is the cyclic gradient defined by $\mathbf{D}P = (D_1 P, \ldots, D_m P)$ and

$$C(P_1, \ldots, P_m, Q_1, \ldots, Q_m)$$
$$= \sum_{k,l=1}^m \left( \mu \otimes \mu[\partial_k P_l \times \partial_l Q_k] + \mu(\partial_l \circ \partial_k V \sharp(P_k, Q_l)) \right) + \sum_{k=1}^m \mu(P_k Q_k).$$

We now give a more explicit formula for $\sigma^2(P, Q)$. We therefore need to
study $C$ and the commutation relations of the cyclic gradient and $\Xi$.

Let us define the following operators on $\mathbb{C}\langle X_1, \ldots, X_m \rangle$:

$$\bar{\Xi}_1 : P \longrightarrow \sum_{k=1}^m \partial_k P \sharp D_k V, \qquad \bar{\Xi}_2 : P \longrightarrow \sum_{i=1}^m (I \otimes \mu) M \circ \partial_i^2 P,$$

where $M(A \otimes B \otimes C) = AC \otimes B$. We also define $\bar{\Xi}_0 = \Sigma^{-1} - \bar{\Xi}_2$ and $\bar{\Xi}$ the operator on $\mathbb{C}_0\langle X_1, \ldots, X_m \rangle$ given by $\bar{\Xi}P = \bar{\Xi}_0 P + \bar{\Xi}_1 P$ if $P \in \mathbb{C}_0\langle X_1, \ldots, X_m \rangle$.
We extend $\bar{\Xi}$ to $\mathbb{C}\langle X_1, \ldots, X_m \rangle$ by setting $\bar{\Xi}1 = 0$. We set, for $i = 0, 1, 2$ or
nothing, $\bar{\bar{\Xi}}_{\mathbf{i}}$ the operator on $\mathbb{C}\langle X_1, \ldots, X_m \rangle^m$ such that $\bar{\bar{\Xi}}_{\mathbf{i}}(P_1, \ldots, P_m) = (\bar{\Xi}_i P_1, \ldots, \bar{\Xi}_i P_m)$.

LEMMA 5.1. *For all $l \in \{1, \ldots, m\}$, for all $P \in \mathbb{C}_0\langle X_1, \ldots, X_m \rangle$, the fol-
lowing equalities hold:*

$$D_l \Sigma^{-1} P = \Sigma^{-1} D_l P + D_l P,$$

$$D_l \Xi_1 P = \bar{\Xi}_1 D_l \Sigma P + \sum_{i=1}^m \partial_i D_l V \sharp D_i \Sigma P,$$

$$D_l \Xi_2 P = \bar{\Xi}_2 D_l \Sigma P.$$

*Besides, let* $\mathrm{Hess}(V) : \mathbb{C}\langle X_1, \ldots, X_m \rangle^m \rightarrow \mathbb{C}\langle X_1, \ldots, X_m \rangle^m$ *be given by*

$$\mathrm{Hess}(V)(v)_l = \sum_{i=1}^m \partial_i D_l V \sharp v_i.$$

*Then, for any* $(P_1, \ldots, P_m) \in \mathbb{C}\langle X_1, \ldots, X_m \rangle^m$, *with $I$ the identity on $\mathbb{C}\langle X_1,$
$\ldots, X_m \rangle^m$, the following relation of commutation relation holds:*

$$\mathbf{D}\Xi = (I + \mathrm{Hess}(V) + \bar{\bar{\Xi}})\mathbf{D}\Sigma.$$

Proof. By linearity, it is sufficient to prove these equalities for a monomial $P = X_{i_1} \cdots X_{i_p}$. Moreover, the projection $\Pi$ onto $\mathbb{C}_0\langle X_1, \ldots, X_m \rangle$ is irrelevant in the definition of the operators $\Xi_i$'s since they are followed by derivatives:

$$D_l \Sigma^{-1} P = p D_l P = (p-1) D_l P + D_l P = \Sigma^{-1} D_l P + D_l P.$$

To prove the second equality, write

$$D_l \Xi_1 P = D_l \sum_{i, \Sigma P = q_1 X_i q_2} q_1 D_i V q_2,$$

then $D_l$ can differentiate $q_1$, $q_2$ or $D_i V$ so that

$$D_l \Xi_1 P = \sum_{i, \Sigma P = r_1 X_l r_2 X_i r_3} r_2 D_i V r_3 r_1 + \sum_{i, \Sigma P = r_1 X_i r_2 X_l r_3} r_3 r_1 D_i V r_2$$

$$+ \sum_{i, \Sigma P = q_1 X_i q_2, D_i V = q_3 X_l q_4} q_4 q_2 q_1 q_3.$$

The sum of the first two terms gives exactly $\bar{\Xi}_1 D_l \Sigma P$ and the last one is

$$\sum_{i, D_i V = q_3 X_l q_4} q_4 D_i P q_3 = \partial_i D_l V \sharp D_i \Sigma P.$$

Note that if $P$ is a monomial,

$$\Xi_2 P = 2 \sum_{i, \Sigma P = q_1 X_i q_2 X_i q_3} \{ \mu[q_1 q_3] q_2 + \mu[q_2] q_1 q_3 \}$$

so that we obtain

$$D_l \Xi_2 P = 2 \sum_{i, \Sigma P = q_1 X_l q_1' X_i q_2 X_i q_3} \mu[q_2] q_1' q_3 q_1$$

$$+ 2 \sum_{i, \Sigma P = q_1 X_i q_2 X_l q_2' X_i q_3} \mu[q_3 q_1] q_2' q_2$$

$$+ 2 \sum_{i, \Sigma P = q_1 X_i q_2 X_i q_3 X_l q_3'} \mu[q_2] q_3' q_1 q_3.$$

Similar algebra shows that

$$\bar{\Xi}_2 D_l \Sigma P = 2 \sum_{i, D_l \Sigma P = q_1 X_i q_2 X_i q_3} \{ \mu(q_2) q_3 q_1 \} = D_l \Xi_2 P.$$

Finally, the last point we only have to sum the previous equalities for $P \in \mathbb{C}_0\langle X_1, \ldots, X_m \rangle$ and all $l \in \{1, \ldots, m\}$,

$$D_l \Xi \Sigma^{-1} P = (D_l + \Sigma^{-1} D_l - \bar{\Xi}_2 D_l + \bar{\Xi}_1 D_l)(P) + \sum_{i=1}^{m} \partial_i D_l V \sharp D_i P$$

$$= [(I + \mathrm{Hess}(V) + \bar{\boldsymbol{\Xi}}) \mathbf{D} P]_l. \qquad \square$$

Thus, we can deduce an expression for $\mathbf{D} \circ \Sigma \Xi^{-1}$.

LEMMA 5.2.   *The operator $\bar{\bar{\Xi}}$ is a symmetric nonnegative operator in $L^2(\mu)$. Let $\cdot^t$ be the involution on $\mathbb{C}\langle X_1, \ldots, X_m \rangle \otimes \mathbb{C}\langle X_1, \ldots, X_m \rangle$ defined by $(A \otimes B)^t = B \otimes A$, then for any $(P, Q) \in \mathbb{C}\langle X_1, \ldots, X_m \rangle^m$,*

$$\mu(P\bar{\bar{\Xi}}Q) = \sum_{k=1}^{m} \mu \otimes \mu(\partial_k P \times [\partial_k Q]^t).$$

$\bar{\bar{\Xi}}$ *is thus nonnegative in $L^2(\mu)^m$ equipped with the scalar product $\langle \mathbf{P}, \mathbf{Q} \rangle = \sum_{i=1}^{m} \mu(P_i Q_i^*)$. $\frac{1-c}{2} I + \operatorname{Hess} V$ is a nonnegative operator in the sense that for every polynomial $P_1, \ldots, P_m$,*

$$\sum_{i=1}^{m} (\operatorname{Hess}(V)P)_i P_i^* \geq -(1-c) \sum_{i=1}^{m} P_i P_i^*.$$

*Thus, $(I + \operatorname{Hess} V + \bar{\bar{\Xi}})$ is symmetric definite positive in $L^2(\mu)^m$ and is invertible. If we consider $\mathbf{D}\Sigma\Xi^{-1}$ as a continuous operator from $\mathbb{C}_0\langle X_1, \ldots, X_m \rangle_A$ into $L^2(\mu)^m$, the following rule of commutation holds:*

$$\mathbf{D}\Sigma\Xi^{-1} = (I + \operatorname{Hess} V + \bar{\bar{\Xi}})^{-1}\mathbf{D}.$$

PROOF.   Here, it is easier to come back to the origin of the problem. The idea is that the operator $\bar{\Xi}$ is a projection of the Laplace operator

$$L = \frac{1}{N} \sum_{k=1}^{m} \sum_{i,j=1}^{N} e^{N\operatorname{Tr}(V+2^{-1}\sum X_l^2)} \partial_{x_{ij}^k} e^{-N\operatorname{Tr}(V+2^{-1}\sum X_l^2)} \partial_{x_{ji}^k}$$

on functions of the matrices. Here, $\partial_{x_{ji}^k}$ is a notation and stands for

$$\tfrac{1}{2}(\partial_{\Re e x_{ji}^k} + \sqrt{-1}\partial_{\Im m x_{ji}^k}).$$

In fact, if we take $P$ a polynomial function,

$$LP = \frac{1}{N}\left[ \sum_{k=1}^{m} \sum_{i,j=1}^{N} N(-D_k V - X_k)_{ji}\, \partial_k P \sharp \Delta_{ji} + \sum_{k=1}^{m} \sum_{i,j=1}^{N} \partial_k \circ \partial_k P \sharp(\Delta_{ij}, \Delta_{ji}) \right]$$

$$= \sum_{k=1}^{m} \partial_k P \sharp(-D_k V - X_k) + \sum_{k=1}^{m} (I \otimes \hat{\mu}^N) M \circ (\partial_k \circ \partial_k) P,$$

with $\Delta_{ij}$ the matrix with null entries except in $(i,j)$ where it is equal to 1. As a consequence, we deduce from the convergence of $\hat{\mu}^N$ toward $\mu$ that, for all polynomials $P, Q$,

$$\lim_{N \to \infty} \int \frac{1}{N} \operatorname{Tr}(QLP)\, d\mu_N^V = -\mu(Q\bar{\Xi}P).$$

But now, by integration by parts, we obtain

$$\int \frac{1}{N} \mathrm{Tr}(QLP) \, d\mu_N^V$$

$$= \int \frac{1}{N^2} \sum_{\alpha,\beta=1}^{N} Q_{\alpha,\beta} (LP)_{\beta,\alpha} \, d\mu_N^V$$

(20)
$$= -\int \frac{1}{N^2} \sum_{k=1}^{m} \sum_{i,j,\alpha,\beta=1}^{N} \partial_{x_k^{ij}} Q_{\alpha,\beta} \, \partial_{x_k^{ji}} P_{\beta,\alpha} \, d\mu_N^V$$

$$= -\int \frac{1}{N^2} \sum_{k=1}^{m} \sum_{i,j,\alpha,\beta=1}^{N} [\partial_k Q \sharp \Delta_{ij}]_{\alpha,\beta} [\partial_k P \sharp \Delta_{ij}]_{\beta,\alpha} \, d\mu_N^V$$

$$= -\sum_{k=1}^{m} \int \hat{\mu}^N \otimes \hat{\mu}^N (\partial_k P \times (\partial_k Q)^t) \, d\mu_N^V,$$

which converges as $N$ goes to infinity toward

$$\sum_{k=1}^{m} \mu \otimes \mu (\partial_k P \times (\partial_k Q)^t) = \mu(Q \bar{\Xi} P).$$

This shows that $\bar{\Xi}$ is symmetric and nonnegative [since if $Q = P^*$, the right-hand side of (20) is clearly nonpositive for all $N$]. Similarly, remark that

$$(\mathrm{Hess}\, V P)_l = \sum_i \partial_i D_l V \sharp P_i.$$

Once estimated at a finite matrix, it is easily seen that

$$\mathrm{Tr}(\partial_i D_l V \sharp P_i P_l^*) = \sum_{\alpha,\beta,\gamma,\delta} (\partial_{x_{\alpha\beta}^i} \partial_{x_{\gamma\delta}^l} \mathrm{Tr} V)(P_i)_{\beta\alpha} (\overline{P_l})_{\delta\gamma}$$

and so the positivity of Hess is deduced at finite $N$ from the convexity of $V$ which, by definition, is the positivity of the Hessian of $\mathrm{Tr}(V)$ in any finite dimension. As a consequence, the operator $I + \mathrm{Hess}(V) + \bar{\bar{\Xi}}$ is invertible on $\mathbb{C}\langle X_1, \ldots, X_m \rangle^m \subset (L^2(\mu))^m$. We then obtain the commutation relation by using the third point of the previous lemma. □

This gives us an explicit formula for $\sigma^2$.

LEMMA 5.3. *For all $P, Q$ in $\mathbb{C}_0\langle X_1, \ldots, X_m \rangle$, for all $1 \le k, l \le m$, the following identities hold:*

$$\mu \otimes \mu[\partial_k D_l P \times \partial_l D_k Q] = \mu \otimes \mu[\partial_k D_l P \times [\partial_k D_l Q]^t],$$

$$C(\mathbf{D}P, \mathbf{D}Q) = \sum_{i=1}^{m} \mu(D_i P[(I + \mathrm{Hess}\, V + \bar{\bar{\Xi}})\mathbf{D}Q]_i),$$

$$\sigma^2(\Xi P, Q) = \sum_{i=1}^m \mu(D_i \Sigma P D_i Q),$$

$$\sigma^2(P, Q) = \sum_{i=1}^m \mu(D_i P (I + \mathrm{Hess}\, V + \bar{\bar{\Xi}})^{-1} D_i Q).$$

PROOF. An elementary computation shows that, for all polynomials $P$,

$$\partial_k D_l P = (\partial_l D_k P)^t.$$

To prove the second equality, recall that

$$C(\mathbf{D}P, \mathbf{D}Q) = \sum_{k,l=1}^m \left( \mu \otimes \mu[\partial_k D_l P \times \partial_l D_k Q] + \mu(\partial_l \circ \partial_k V \sharp(D_k P, D_l Q)) \right)$$

$$+ \sum_{k=1}^m \mu(D_k P D_k Q).$$

The third term can be directly written $\sum_{i=1}^m \mu(D_i P[\mathbf{D}Q]_i)$. For the second term, we use the first equality and Lemma 5.2:

$$\sum_{k,l=1}^m \mu \otimes \mu[\partial_k D_l P \times \partial_l D_k Q] = \sum_{i=1}^m \mu(D_i P \bar{\Xi} D_i Q).$$

Finally, we only need to check if the two terms in the second derivative of $V$ coincide, but this is clear by the trace property:

$$\sum_{k,l=1}^m \mu(\partial_l \circ \partial_k V \sharp(D_k P, D_l Q)) = \sum_{i,j=1}^m \mu(D_i P \partial_j D_i V \sharp D_j Q).$$

For the last points we only have to use the commutation rule of Lemma 5.2 and the previous point:

$$\sigma^2(\Xi P, Q) = C(\mathbf{D}\Sigma P, \mathbf{D}\Sigma \Xi^{-1} Q)$$

$$= \sum_{i=1}^m \mu(D_i \Sigma P[(I + \mathrm{Hess}\, V + \bar{\Xi})\mathbf{D}\Sigma \Xi^{-1} Q]_i)$$

$$= \sum_{i=1}^m \mu(D_i \Sigma P D_i Q).$$

The last point is proved with the same technique. $\square$

5.2. *Combinatorial interpretation.* It was shown in [14] that for small **t**'s the limit measure $\mu$ has a combinatorial interpretation. More precisely, let $V = \sum_i t_i q_i$ with some monomials $q_i$. Note that in order to have a self-adjoint potential, in the decomposition in monomials, the coefficient of a monomial must be the complex conjugate of the coefficient of its adjoint.

We define a set of colors as the set $\{1, \ldots, m\}$ and associate to each monomial $q = X_{i_1} \cdots X_{i_p}$ a star (i.e., a vertex with some half-edges pointing out of it) of $p$ half-edges which are in the clockwise order respectively of color $i_1$, $i_2, \ldots, i_p$. Besides, we distinguish the first half-edge so that we clearly obtain a bijection between monomials and stars. We will say that the star is of type $q$ if it comes from a monomial $q$ in that way. Note that a star can equivalently be represented by an annulus with ordered colored dots and a distinguished dot.

Given a set of such stars embedded in the sphere, we can construct some graphs among them simply by gluing pairwise different half-edges of the same color and such that the resulting edges do not cross each other. We call a graph obtained in this way a planar graph. Two planar graphs are said to be equivalent if there is a homeomorphism of the sphere which fix each star and take the first graph on the second. A map is a class of equivalence of connected planar graphs for the relation of homomorphism. We now define

$$\mathcal{M}_{k_1,\ldots,k_n}(P) = \sharp \left\{ \begin{array}{c} \text{maps with } k_i \text{ stars of type } q_i \\ \text{and one of type } P \end{array} \right\}$$

and

$$\mathcal{M}_{k_1,\ldots,k_n}(P,Q) = \sharp \left\{ \begin{array}{c} \text{maps with } k_i \text{ stars of type } q_i \\ \text{one of type } P \text{ and one of type } Q \end{array} \right\}.$$

These quantities are only defined for $P$ and $Q$ monomials, but we immediately extend them by linearity to arbitrary polynomials $P$ and $Q$. By convention, the star associated to the monomial 1 is empty so that $\mathcal{M}_{k_1,\ldots,k_n}(P,1) = 0$.

In [14], Section 3.2 there is the following relation between the limit measure and the enumeration of planar graphs.

THEOREM 5.4. *There exists $\eta > 0$ such that, for $\mathbf{t} \in B_\eta$, for all polynomial $P$,*

$$\mu(P) = \sum_{k_1,\ldots,k_n} \prod_{i=1}^{n} \frac{(-t_i)^{k_i}}{k_i!} \mathcal{M}_{k_1,\ldots,k_n}(P).$$

We now prove that there is a similar link between the variance $\sigma^2(P)$ which appears in our central limit theorem and the generating function of the $\mathcal{M}_{k_1,\ldots,k_n}(P,Q)$. We define

$$\mathcal{M}(P,Q) = \sum_{k_1,\ldots,k_n} \prod_{i=1}^{n} \frac{(-t_i)^{k_i}}{k_i!} \mathcal{M}_{k_1,\ldots,k_n}(P,Q).$$

We shall prove that $\sigma^2(P,Q)$ and $\mathcal{M}(P,Q)$ satisfy the same kind of induction relation.

PROPOSITION 5.5.   *For all monomials $P,Q$ and all $k$,*

$$\mathcal{M}_{k_1,\ldots,k_n}(X_k P,Q)$$
$$= \sum_{0 \le p_i \le k_i} \sum_{P=RX_kS} \prod_i C_{k_i}^{p_i} \mathcal{M}_{p_1,\ldots,p_n}(R,Q) \mathcal{M}_{k_1-p_1,\ldots,k_n-p_n}(S)$$
$$+ \sum_{0 \le p_i \le k_i} \sum_{P=RX_kS} \prod_i C_{k_i}^{p_i} \mathcal{M}_{p_1,\ldots,p_n}(S,Q) \mathcal{M}_{k_1-p_1,\ldots,k_n-p_n}(R)$$
$$+ \sum_{0 \le j \le n} k_j \mathcal{M}_{k_1,\ldots,k_j-1,\ldots,k_n}(D_k V P,Q) + \mathcal{M}_{k_1,\ldots,k_n}(D_k Q P)$$

*and*

$$(21) \quad \mathcal{M}(X_k P,Q) = \mathcal{M}((I \otimes \mu + \mu \otimes I)\partial_k P) - \mathcal{M}(D_k V P,Q) + \mu(D_k Q P).$$

*Besides, there exists $\eta > 0$ so that, there exists $R < +\infty$ such that for all monomials $P$ and $Q$, all $\mathbf{t} \in B(0,\eta)$,*

$$|\mathcal{M}(P,Q)| \le R^{\deg P + \deg Q}.$$

PROOF.   The proof is very close to that given of Theorem 2.2 in [14] which explains the decomposition of planar maps with one root. We look at the first half-edge with color $k$ corresponding to $X_k$ in $X_k P$:

1. The first possibility is that the half-edge is glued to another half-edge of $P = R X_k S$. It cuts $P$ in two monomials $R$ and $S$ and it occurs for all decomposition of $P$ into $P = R X_k S$ which is exactly what does $D$. Then either the component $R$ is linked to $Q$ and to $p_i$ stars of type $q_i$ for each $i$, this leads to

$$\prod_i C_{k_i}^{p_i} \mathcal{M}_{p_1,\ldots,p_n}(R,Q) \mathcal{M}_{k_1-p_1,\ldots,k_n-p_n}(S)$$

   possibilities, or we are in the symmetric case with $S$ linked to $Q$ in place of $R$.

2. The second case occurs when the half-edge is glued to a star of type $q_j$ for a given $j$, then first we have to choose between the $k_j$ vertices of this type, then we contract the edges arising from this gluing to form a star of type $D_i q_j P_1$; there are

$$k_j \mathcal{M}_{k_1,\ldots,k_j-1,\ldots,k_n}(D_k q_j P,Q)$$

   choices.

3. The last case is that the half-edge can be glued with the star associated to $Q = RX_iS$. We contract this half-edge and obtain a star of type $D_kQP$. This leads to

$$\mathcal{M}_{k_1,\dots,k_n}(D_kQP)$$

possibilities.

We can now sum on the $k$'s to obtain the relation on $\mathcal{M}$.

Finally, to show the last point of the proposition, we only have to prove that there exists $A > 0, B > 0$ such that, for all $k$'s, for all monomials $P$ and $Q$,

$$\frac{\mathcal{M}_{k_1,\dots,k_n}(P,Q)}{\prod_i k_i!} \le A^{\sum_i k_i} B^{\deg P + \deg Q}.$$

This follows easily by induction over the degree of $P$ with the previous relation on the $\mathcal{M}$ since we have proved such a control for $\mathcal{M}_{k_1,\dots,k_n}(Q)$ in [14]. $\quad\square$

We can now relate the variance and the generating function for the enumeration of planar maps with two prescribed vertices.

THEOREM 5.6. *Assume* $(\mathbf{H}')$ *with $\eta$ small enough. Then, for all polynomials $P, Q$,*

$$\sigma^2(P,Q) = \mathcal{M}(P,Q).$$

PROOF. First we transform the relation on $\mathcal{M}$. We use (21) with $P = D_k\Sigma R$ to deduce

$$\mathcal{M}(\Xi R, Q) = \sum_k \mu(D_kQD_k\Sigma R).$$

Let us define $\Delta = \sigma^2 - \mathcal{M}$. Then according to (5.2) and the previous property, $\Delta$ is compactly supported and for all polynomials $P$ and $Q$,

$$\Delta(\Xi P, Q) = 0.$$

Moreover, with $\mathcal{M}(1,Q) = 0 = \sigma^2(1,Q)$,

$$\Delta(1,Q) = 0.$$

To conclude, we have to invert one more time the operator $\Xi$. For a polynomial $P$, we take, as in the proof of the central limit theorem, a sequence of polynomial $S_n$ which goes to $S = \Xi^{-1}P$ in $\mathbb{C}_0\langle X_1,\dots,X_m\rangle_A$. Then, write

$$\Delta(P,Q) = \Delta(\Xi(S_n + S - S_n), Q) = \Delta(\Xi(S - S_n), Q).$$

But by continuity of $\Xi$, $\Xi(S - S_n)$ goes to 0 for the norm $\|\cdot\|_A$. We can always assume $A \ge R$ if $\eta$ is small enough. Moreover, because $\Delta$ is compactly supported, $\Delta$ is continuous for $\|\cdot\|_A$, and so $\Delta(\Xi(S - S_n), Q)$ goes to zero when $n$ goes to $+\infty$. This proves the theorem. $\quad\square$

**6. Second order correction to the free energy.** We now deduce from the central limit theorem the precise asymptotics of $N\overline{\delta}^N(P)$ and then compute the second order correction to the free energy.

Let $\phi_0$ and $\phi$ be the linear forms on $\mathbb{C}_0\langle X_1,\ldots,X_m\rangle$ which are given, if $P$ is a monomial by

$$(22) \qquad \phi_0(P) = \sum_{i=1}^m \sum_{P=P_1 X_i P_2 X_i P_3} \sigma^2(P_3 P_1, P_2)$$

and $\phi = \phi_0 \circ \Sigma$.

PROPOSITION 6.1. *Assume* (**H′**). *Then, for any* $P$ *in* $\mathbb{C}_0\langle X_1,\ldots,X_m\rangle$,

$$\lim_{N\to\infty} N\overline{\delta}^N(P) = \phi(\Xi^{-1}\Pi(P)).$$

PROOF. Again, we base our proof on the finite dimensional Schwinger–Dyson equation (13) which, after centering, and since we can always assume that $P \in \mathbb{C}_0\langle X_1,\ldots,X_m\rangle$, reads for $i \in \{1,\ldots,m\}$,

$$N^2\mu_V^N((\hat{\mu}^N - \mu)[(X_i + D_i V)P - (I\otimes\mu + \mu\otimes I)\partial_i P]) = \mu_V^N(\hat{\delta}^N\otimes\hat{\delta}^N(\partial_i P)).$$

Taking $P = D_i \Sigma P$ and summing over $i \in \{1,\ldots,m\}$, we thus have

$$(23) \qquad N^2\mu_V^N((\hat{\mu}^N - \mu)(\Xi P)) = \mu_V^N\left(\hat{\delta}^N\otimes\hat{\delta}^N\left(\sum_{i=1}^m \partial_i \circ D_i \Sigma P\right)\right).$$

By Corollary 4.8 and Lemma 5.1, we see that

$$\lim_{N\to\infty} \mu_V^N\left(\hat{\delta}^N\otimes\hat{\delta}^N\left(\sum_{i=1}^m \partial_i \circ D_i \Sigma P\right)\right) = \phi(P),$$

which gives the asymptotics of $N\overline{\delta}^N(\Xi P)$ for all $P$.

To generalize the result to arbitrary $P$, we proceed as in the proof of the full central limit theorem. We take a sequence of polynomials $Q_n$ which goes to $Q = \Xi^{-1}P$ when $n$ goes to $\infty$ for the norm $\|\cdot\|_A$. We denote $R_n = P - \Xi Q_n = \Xi(Q - Q_n)$. Note that as $P$ and $Q_n$ are polynomials, then $R_n$ is also a polynomial. Then we write

$$N\overline{\delta}^N(P) = N\overline{\delta}^N(\Xi Q_n) + N\overline{\delta}^N(R_n).$$

According to Property 3.1, for any monomial $P$ of degree less than $\varepsilon N^{2/3}$,

$$|N\overline{\delta}^N(P)| \le C^{\deg(P)}.$$

So if we take the limit in $N$, for any monomial $P$,

$$\limsup_N |N\overline{\delta}^N(P)| \le C^{\deg(P)}$$

and if $P$ is a polynomial,

$$\limsup_N |N\overline{\delta}^N(P)| \leq \|P\|_C \leq \|P\|_A.$$

The last inequality comes from the hypothesis $(\mathbf{H}')$ which require $C < A$.

We now fix $n$ and let $N$ go to infinity,

$$\limsup_N |N\overline{\delta}^N(P - \Xi Q_n)| \leq \limsup_N |N\overline{\delta}^N(R_n)| \leq \|R_n\|_A.$$

If we now let $n$ go to infinity, the right-hand side term vanishes and we are left with

$$\lim_N N\overline{\delta}^N(P) = \lim_n \lim_N N\overline{\delta}^N(Q_n) = \lim_n \phi(Q_n).$$

It is now sufficient to show that $\phi$ is continuous for the norm $\|\cdot\|_A$. But it can be checked easily that $P \to \sum_{i=1}^m \partial_i \circ D_i P$ is continuous from $\mathbb{C}_0\langle X_1,\ldots,X_m\rangle_A$ to $\mathbb{C}_0\langle X_1,\ldots,X_m\rangle_{A-1}$ and $\sigma^2$ is continuous for $\|\cdot\|_{A-1}$ due to the technical hypothesis in $(\mathbf{H}')$. This proves that $\phi$ is continuous and then can be extended on $\mathbb{C}_0\langle X_1,\ldots,X_m\rangle_A$. Thus,

$$\lim_N N\overline{\delta}^N(P) = \lim_n \phi(Q_n) = \phi(Q). \qquad \square$$

This result allows us to estimate the first order correction to the free energy.

THEOREM 6.2. *Assume* $(\mathbf{H}')$, *then the following asymptotics hold:*

$$\log Z_{V_{\mathbf{t}}}^N = N^2 F^0(V_{\mathbf{t}}) + F^1(V_{\mathbf{t}}) + o(1),$$

*with*

$$F^0(V_{\mathbf{t}}) = -\int_0^1 \mu_{\alpha\mathbf{t}}\left(\sum_{i=1}^n t_i q_i\right) d\alpha$$

*and*

$$F^1(V_{\mathbf{t}}) = -\int_0^1 \phi_{\alpha\mathbf{t}}\left(\Xi_{\alpha\mathbf{t}}^{-1} \sum_{i=1}^n t_i q_i\right) ds,$$

*with* $\Xi_{\alpha\mathbf{t}}$ *(resp.* $\phi_{\alpha\mathbf{t}}$*) the operator* $\Xi$ *(resp. the linear form* $\phi$*) corresponding to the potential* $V_{\alpha\mathbf{t}} = \alpha V_{\mathbf{t}}$ *with parameters* $\alpha\mathbf{t}$.

PROOF. Remark that, for $i \in \{1,\ldots,n\}$,

$$\partial_\alpha \log Z_{\alpha V_{\mathbf{t}}}^N = -N^2 \mu_{\alpha V_{\mathbf{t}}}^N\left(\hat{\mu}^N\left(\sum_{i=1}^n t_i q_i\right)\right)$$

so that we can write

$$(24) \qquad \log Z_{V_{\mathbf{t}}}^N = N^2 F^0(V_{\mathbf{t}}) - \int_0^1 \left[ N \bar\delta_{\alpha \mathbf{t}}^N \left( \sum t_i q_i \right) \right] d\alpha.$$

Since for all $\alpha \in [0, 1]$, $V_{\alpha \mathbf{t}} = \alpha V_{\mathbf{t}}$ is $c \wedge 1$-convex if $V_{\mathbf{t}}$ is $c$-convex, Proposition 6.1 and (24) finish the proof of the theorem since, by Proposition 3.1, all the $N \bar\delta_{\alpha \mathbf{t}}^N (q_i)$ can be bounded independently of $N$, $\alpha \in [0, 1]$ and $t \in B_{\eta,c}$ so that the dominated convergence theorem applies. $\quad \square$

As for the combinatorial interpretation of the variance, we relate $F^1(V_{\mathbf{t}})$ to a generating function of maps. This time, we will consider maps on a torus instead of a sphere. Such maps are said to be of genus 1. We define

$$\mathcal{M}_{k_1,\dots,k_n}^1(P) = \sharp \left\{ \begin{array}{c} \text{maps of genus 1 with } k_i \text{ stars of type } q_i \\ \text{and one of type } P \end{array} \right\}$$

and

$$\mathcal{M}_{k_1,\dots,k_n}^1 = \sharp \{ \text{maps with } k_i \text{ stars of type } q_i \}.$$

We also define the generating function

$$\mathcal{M}^1(P) = \sum_{k_1,\dots,k_n} \prod_{i=1}^n \frac{(-t_i)^{k_i}}{k_i!} \mathcal{M}_{k_1,\dots,k_n}^1(P).$$

If $P$ is a monomial, we will denote $\mathcal{M}(\partial_i P)$ for $\sum_{P=RX_iS} \mathcal{M}(R,S)$ and we extend this notation to all polynomials by linearity.

PROPOSITION 6.3.  *For all monomials $P$ and all $k$,*

$$\mathcal{M}_{k_1,\dots,k_n}^1(X_k P)$$
$$= \sum_{0 \le p_i \le k_i} \sum_{P=RX_kS} \prod_i C_{k_i}^{p_i} \mathcal{M}_{p_1,\dots,p_n}^1(R) \mathcal{M}_{k_1-p_1,\dots,k_n-p_n}(S)$$
$$\quad + \sum_{0 \le p_i \le k_i} \sum_{P=RX_kS} \prod_i C_{k_i}^{p_i} \mathcal{M}_{p_1,\dots,p_n}(R) \mathcal{M}_{k_1-p_1,\dots,k_n-p_n}^1(S)$$
$$\quad + \sum_{0 \le j \le n} k_j \mathcal{M}_{k_1,\dots,k_j-1,\dots,k_n}^1(D_k V P, Q) + \sum_{P=RX_kS} \mathcal{M}_{k_1,\dots,k_n}(R,S)$$

*and*

$$(25) \qquad \mathcal{M}^1(X_k P) = \mathcal{M}^1((I \otimes \mu + \mu \otimes I) \partial_k P) - \mathcal{M}^1(D_k V P) + \mathcal{M}(\partial_k P).$$

*Besides, for $\eta$ small enough, there exists $R < +\infty$ such that, for all monomials $P$, all $\mathbf{t} \in B(0,\eta)$,*

$$|\mathcal{M}^1(P)| \le R^{\deg P}.$$

PROOF. We proceed as we did for the combinatorial interpretation of the variance. We look at the first half-edge corresponding to $X_k$, then two cases may occur.

1. The first possibility is that the half-edge is glued to another half-edge of $P = RX_kS$. It forms a loop starting from $P$. There are two cases:

(a) The loop can be retractable. It cuts $P$ in two monomials $R$ and $S$ and it occurs for all decomposition of $P$ into $P = RX_kS$ which is exactly what does $D$. Then either the component $R$ or the component $S$ is of genus 1 and the other component is planar. It produces either

$$\prod_i C_{k_i}^{p_i} \mathcal{M}^1_{p_1,\ldots,p_n}(R)\mathcal{M}_{k_1-p_1,\ldots,k_n-p_n}(S)$$

possibilities or the symmetric formula (where we exchange $R$ and $S$).

(b) The loop can also be nontrivial in the fundamental group of the surface. Then the surface is cut in two. We are left with a planar surface with two fixed stars $R$ and $S$. This gives

$$\mathcal{M}_{k_1,\ldots,k_n}(R,S)$$

possibilities.

2. The second possibility occurs when the half-edge is glued to a half-edge of a star of type $q_j$ for a given $j$, then first we have to choose between the $k_j$ stars of this type, then we contract the edges arising from this gluing to form a star of type $D_iq_jP_1$; this creates

$$k_j\mathcal{M}^1_{k_1,\ldots,k_j-1,\ldots,k_n}(D_kq_jP,Q)$$

possibilities.

We can now sum on the $k$'s to obtain the relation on $\mathcal{M}^1$.

Finally, to show that $\mathcal{M}^1$ is compactly supported, we only have to prove that there exists $A > 0, B > 0$ such that, for all $k$'s, for all monomials $P$,

$$\frac{\mathcal{M}^1_{k_1,\ldots,k_n}(P)}{\prod_i k_i!} \leq A^{\sum_i k_i} B^{\deg P}.$$

Another time this follows easily by induction with the previous relation on the $\mathcal{M}^1(P)$'s. $\square$

We now give the combinatorial interpretation for the first order correction to the free energy.

PROPOSITION 6.4. Assume $(\mathbf{H'})$. There exists $\eta > 0$ small enough so that, for $\mathbf{t} \in B_{\eta,c}$, for all nonconstant monomial $P$,

$$\phi(\Xi^{-1}P) = \mathcal{M}^1(P)$$

*and*

$$F^1 = \sum_{k_1,\ldots,k_n \in \mathbb{N}^n - \{0\}} \prod_{i=1}^n \frac{(-t_i)^{k_i}}{k_i!} \mathcal{M}^1_{k_1,\ldots,k_n}.$$

PROOF.   We use the previous property with $P = D_k \Sigma P$ and we sum on $k$:

$$\mathcal{M}^1(\Xi P) = \mathcal{M}\left(\sum_k \partial_k D_k \Sigma P\right) = \sum_k \sigma^2(\partial_k D_k \Sigma P) = \phi(P),$$

where we have used the combinatorial interpretation of the variance (Theorem 5.6). As $\mathcal{M}^1$ and $\phi$ are continuous for $\|\cdot\|_A$ when $\eta$ is small enough, we can apply this to $\Xi^{-1}P$ and conclude.

Finally, for $\eta$ sufficiently small, the series is absolutely convergent so that we can invert the integral and the sum to obtain

$$\begin{aligned}
F^1(V_{\mathbf{t}}) &= -\int_0^1 \mathcal{M}^1_{\alpha t_1,\ldots,\alpha t_n}\left(\sum_j t_j q_j\right) d\alpha \\
&= \int_0^1 \sum_{k_1,\ldots,k_n} \sum_j \prod_i \frac{(-\alpha t_i)^{k_i}}{k_i!}(-t_j)\mathcal{M}^1_{k_1,\ldots,k_n}(q_j)\, d\alpha \\
&= \sum_{k_1,\ldots,k_n} \frac{1}{k_1+\cdots+k_n+1} \sum_j \prod_i \frac{(-t_i)^{k_i}}{k_i!}(-t_j)\mathcal{M}^1_{k_1,\ldots,k_j+1,\ldots,k_n} \\
&= \sum_{k_1,\ldots,k_n} \prod_i \frac{(-t_i)^{k_i}}{k_i!}\mathcal{M}^1_{k_1,\ldots,k_j,\ldots,k_n}. \qquad \square
\end{aligned}$$

**7. Diverging integrals.**   Physicists often use matrix models in more general settings. We would like to study the case of a potential $V$ for which the integral $Z_V^N$ is not convergent. For example, one may wonder if we can obtain the generating function for planar triangulation. The issue is that for $V = tX^3$, $Z_V^N$ is infinite. The idea to give a meaning to this integral is to add a cut-off; we define, for $L > 0$,

$$Z_{V,L}^N = \int_{\mathcal{H}_N(\mathbb{C})^m, \lambda_{\max}(\mathbf{A})<L} e^{-N\mathrm{Tr}(V(A_1,\ldots,A_m))}\, d\mu^N(A_1,\ldots,A_m).$$

This allows us to define the probability measure

$$\begin{aligned}
\mu_{V,L}^N(dA_1,\ldots,dA_m) \\
= \frac{\mathbb{1}_{\lambda_{\max}(\mathbf{A})<L}}{Z_{V,L}^N} e^{-N\mathrm{Tr}(V(A_1,\ldots,A_m))}\, d\mu^N(A_1,\ldots,A_m).
\end{aligned}$$

In [14], we show that, for all $L > L_0$ for a well chosen $L_0$, there exists $\eta > 0$ such that for $|\mathbf{t}| < \eta$, $\hat{\mu}^N$ goes almost surely toward the unique solution to Schwinger–Dyson's equation (2). This shows that the cut-off does not perturb too much the model since the limit does not depend on the choice of the cut-off $L$ and keeps the same interpretation than in case of convex potentials. The aim of this section is to show that we can also extend the central limit theorem to this setting. The key idea is to see this potential as a convex potential. We bound the Hessian of

$$(26) \qquad \varphi_{V_\mathbf{t}}^N : (A_k(ij)) \in (\mathbb{R}^{N^2})^m \cap \{\lambda_{\max}(\mathbf{A}) \le L\} \to \mathrm{Tr}(V(A_1, \dots, A_m))$$

uniformly in $N$:

$$\mathrm{Hess}\,\varphi_{V_\mathbf{t}}^N(A, A) = \sum_{i=1}^n t_i \sum_{q_i = RXSXT} \mathrm{Tr}(RASAT).$$

Now, using Hölder's inequality,

$$|\mathrm{Tr}(RASAT)| = |\mathrm{Tr}(TRASA)| \le \sqrt{\mathrm{Tr}((TR)A^*A(TR)^*)}\sqrt{\mathrm{Tr}(SA^*AS^*)}$$
$$\le \|TR\|\|S\|\mathrm{Tr}(AA^*),$$

which implies that, for $\{\lambda_{\max}(\mathbf{A}) \le L\}$,

$$\|\mathrm{Hess}\,\varphi_{V_\mathbf{t}}^N\| \le C|\mathbf{t}|$$

and $C$ depends only on $L$. Therefore, we can find $\varepsilon > 0$ such that if $\mathbf{t} \in B(0, \varepsilon) \cap \{\mathbf{t} | V_\mathbf{t} = V_\mathbf{t}^*\}$, for all $N$, $\varphi_{V_\mathbf{t}}^N + \frac{1}{4}\sum_{i=1}^n \mathrm{Tr}(X_i^2)$ is convex on $\{\lambda_{\max}(\mathbf{A}) \le L\}$.

Thus, $\tilde{V}_\mathbf{t}(\mathbf{A}) = V_\mathbf{t}(\mathbf{A}) + \infty \mathbb{1}_{\lambda_{\max}(\mathbf{A}) > L}$ is a convex potential and

$$\mathbb{1}_{\lambda_{\max}(\mathbf{A}) \le L} e^{-N\mathrm{Tr}(V_\mathbf{t}(\mathbf{A}))} = e^{-N\mathrm{Tr}(\tilde{V}(\mathbf{A}))}$$

is log-concave so that most of the step we proved so far can be generalized to this case. Indeed, the Brascamp–Lieb and concentration inequalities do not require smoothness for the potential $V$. In fact, we could have included this case in all of the previous proofs but they would have been less readable. We will only sketch the proof in this generalized case and highlight the main differences with the convex case.

First, we must control the rate of convergence of the measure to its limit. The important fact is that up to the choice of $\mathbf{t}$ we can obtain bounds independent of $L$.

PROPOSITION 7.1. *There exist nonnegative constants $L_0, M_0, C, \alpha$ such that, for $L > L_0$, we can find $\eta > 0$ such that, for $|\mathbf{t}| < \eta$:*

1.
$$\mu(X_i^{2n}) \leq \limsup_N \overline{\mu}^N(X_i^{2n}) \leq C^{2n}.$$

2. *For all $M > M_0$*
$$\mu_V^N(\lambda_{\max}^N(\mathbf{A}) > M) \leq e^{-\alpha M N}.$$

3. *There exists a finite constant $\varepsilon_{P,M}^N$ such that, for any $\varepsilon > 0$,*
$$\mu_V^N(\{|\hat{\underline{\hat{\varrho}}}^N(P) - m_{P,M}^N| \geq \varepsilon + \varepsilon_{P,M}^N\} \cap \Lambda_M^N) \leq 2e^{-c\varepsilon^2/(2\|P\|_{\mathcal{L}}^M)}$$

   *and if $P$ is a monomial of degree $d$, $\varepsilon_{P,M} \leq NCdM^d e^{-\alpha M N}$.*

PROOF.  Since $e^{-N\operatorname{Tr}(\tilde{V}(\mathbf{A}))}$ is log-concave, we can still use the Brascamp–Lieb inequalities. The only point to check is that we can still find a lower bound for $Z_{V,L}^N$, but this was already done in [14] using Jensen's inequality:

$$Z_{V_{\mathbf{t}}}^{N,L} = \int_{\lambda_{\max}(\mathbf{A}) \leq L} e^{-N\operatorname{Tr}(V_{\mathbf{t}}(\mathbf{A}))} \prod d\mu_N(A_i)$$

$$\geq \mu^N(\lambda_{\max}(\mathbf{A}) \leq L) \exp\left(-N \int_{\lambda_{\max}(\mathbf{A}) \leq L} \operatorname{Tr}(V_{\mathbf{t}}(\mathbf{A})) \frac{\prod d\mu_N(A_i)}{\mu^N(\lambda_{\max}(\mathbf{A}) \leq L)}\right).$$

The biggest eigenvalue goes almost surely to 2 and

$$\left|\int_{\lambda_{\max}(\mathbf{A}) \leq L} \frac{1}{N} \operatorname{Tr}(V_{\mathbf{t}}(\mathbf{A})) \prod d\mu_N(A_i)\right|$$

is bounded by $\mu^N(V_{\mathbf{t}}V_{\mathbf{t}}^*)^{1/2}$ which goes to $\sigma^m(V_{\mathbf{t}}V_{\mathbf{t}}^*)^{1/2} < +\infty$ according to [23]. Thus, if $L > 2$, $Z_{V_{\mathbf{t}}}^{N,L} \geq e^{-dN^2}$ for a finite constant $d$. Thus, we can prove the property as in Section 2. The proof of the two last points do not differ from the convex case.  □

The idea, once we have an a priori control on the radius $C$ of the support independently of $L$, is that we can use it to approximate any polynomial by a compactly supported function with support in $[-L, L]$. We choose $L > L_0 = \max(M_0, C)$ and define for $L_0 < R < L$, $\phi_R$ the piecewise affine function such that, for $|x| < R$, $\phi_R(x) = x$ and $\phi_R$ has a compact support strictly inside $[-L, L]$. Then we can approximate any polynomial $P(\mathbf{X})$ by $h_R = P(\phi_R(X_1), \ldots, \phi_R(X_m))$. The main improvement in the replacement of $P$ by $h_R$ is that $h_R$ satisfies the finite Schwinger–Dyson's equation (13).

PROPOSITION 7.2.  *If $L$ is bigger than some $L_0 > 0$, and $\varepsilon > 0$, there exist $C, \eta, M_0$ such that, for $M > M_0$, $|\mathbf{t}| < \eta$ for all polynomial $P$ of degree $d < \varepsilon N^{2/3}$,*
$$|\overline{\delta}^N(P)| \leq C\frac{\|P\|_M}{N}.$$

Proof. In order to prove the analogue in the convex case (Property 3.1), we use the finite Schwinger-Dyson's equation which is not always satisfied in this case. In fact, it is only satisfied for compactly supported function $h$ with support in $[-L, L]$, since for such $h$ we can make the infinitesimal change of variable. For a polynomial $P$,

$$\mu_V^N(\hat{\mu}^N[(X_i + D_iV)P]) - \mu_V^N(\hat{\mu}^N \otimes \hat{\mu}^N(\partial_i P))$$
$$= \mu_V^N(\hat{\mu}^N[(X_i + D_iV)(P - h_R)]) - \mu_V^N(\hat{\mu}^N \otimes \hat{\mu}^N(\partial_i(P - h_R))).$$

Therefore, since $\mu$ satisfies the Schwinger–Dyson equation, we get that, for all polynomial $P$,

$$(27) \qquad \overline{\delta}^N(X_iP) = -\overline{\delta}^N(D_iVP) + \overline{\delta}^N \otimes \overline{\mu}^N(\partial_i P) + \mu \otimes \overline{\delta}^N(\partial_i P) + r(N, P),$$

with

$$r(N, P) := N^{-1}\mu_V^N(\underline{\hat{\delta}}^N \otimes \underline{\hat{\delta}}^N(\partial_i P))$$
$$+ N(\mu_V^N(\hat{\mu}^N[(X_i + D_iV)(P - h_R)])$$
$$- \mu_V^N(\hat{\mu}^N \otimes \hat{\mu}^N(\partial_i(P - h_R)))).$$

Thus, the only difference with the convex case is the term $N(\mu_V^N(\hat{\mu}^N[(X_i + D_iV)(P - h_R)]) - \mu_V^N(\hat{\mu}^N \otimes \hat{\mu}^N(\partial_i(P - h_R))))$ but, since on $\Lambda_M^N$, $P(\mathbf{A}) = h_R(\mathbf{A})$ and $R > M$, this term decreases exponentially fast and this allows to finish the proof exactly as in the proof of Proposition 3.1. □

Since the main tools are available, we next turn to the proof of the central limit theorem. Here we have to be careful since the technique of the "infinitesimal change of variable" is no longer true in its full generality. But it still holds if we restrict ourseves to compactly supported functional, thus, we immediately obtain a weaker version of Lemma 4.1:

LEMMA 7.3. If $L$ is bigger than some $L_0 > 0$, there exists $\eta$ such that, for $|\mathbf{t}| < \eta$ if $h_1, \ldots, h_m$ are compactly supported with support in $]-L, L[$ and self-adjoint, the random variable

$$Y_N(h_1, \ldots, h_m) = N \sum_{k=1}^m \{\hat{\mu}^N \otimes \hat{\mu}^N(\partial_k h_k) - \hat{\mu}^N[(X_k + D_kV)h_k]\}$$

converges in law toward a real centered Gaussian variable with variance

$$C(h_1, \ldots, h_m) = \sum_{k,l=1}^m (\mu \otimes \mu[\partial_k h_l \times \partial_l h_k] + \mu(\partial_l \circ \partial_k V\sharp(h_k, h_l))) + \sum_{k=1}^m \mu(h_k^2).$$

The last step is to show that even if we do not have the result for all Stieljes functions, it is sufficient to approach polynomials by compactly supported function with support inside $]-L, L[$. We will again use the fact that the limit measure has a support bounded independently of $L$. Using this idea, we prove a result similar to Lemma 4.2.

LEMMA 7.4.  *If $L$ is bigger than some $L_0 > 0$, there exists $\eta$ such that for $|\mathbf{t}| < \eta$, if $P_1, \ldots, P_m$ are in $\mathbb{C}\langle X_1, \ldots, X_m \rangle_{sa}$, then the variable*

$$Y_N(P_1, \ldots, P_m) = N \sum_{k=1}^{m} [\hat{\mu}^N \otimes \hat{\mu}^N (\partial_k P_k) - \hat{\mu}^N [(X_k + D_k V) P_k]]$$

*converges in law toward a real centered Gaussian variable with variance*

$$C(P_1, \ldots, P_m) = \sum_{k,l=1}^{m} \left( \mu \otimes \mu [\partial_k P_l \times \partial_l P_k] + \mu(\partial_l \circ \partial_k V \sharp(P_k, P_l)) \right) + \sum_{k=1}^{m} \mu(P_k^2).$$

PROOF.  First choose $L > L_0 = \max(M_0, C)$ and for $L_0 < R < L$, approximate the polynomials $P_i(\mathbf{X})$ by $h_R^i = P_i(\phi_R(X_1), \ldots, \phi_R(X_m))$. Then since $C$ bounds the support of $\mu$, observe that $C(P_1, \ldots, P_m) = C(h_R^1, \ldots, h_R^m)$ and we only have to prove that

$$Y_N(P_1, \ldots, P_m) - Y_N(h_R^1, \ldots, h_R^m)$$

goes in law to 0 when $N$ goes to infinity. But, we have the inequality

$$P(|Y_N(P_1, \ldots, P_m) - Y_N(h_R^1, \ldots, h_R^m)| > \varepsilon) \le P(\lambda_{\max}(\mathbf{A}) > R)$$

and the right-hand side goes exponentially fast to 0.  □

The other results can be proved as in the convex case with only minor modifications. Following the same way than in the convex case, this allows us to prove the theorem:

THEOREM 7.5.  *If $L$ is bigger than some $L_0 > 0$, there exists $\eta$ such that, for $|\mathbf{t}| < \eta$, for all $P$ in $\mathbb{C}_0\langle X_1, \ldots, X_m \rangle$, $\hat{\delta}^N(P)$ converges in law to a Gaussian variable with variance*

$$\sigma^2(P) := \mathcal{C}(\Xi^{-1} P) = C(D_1 \Sigma \Xi^{-1} P, \ldots, D_m \Sigma \Xi^{-1} P).$$

Besides, the convergence in moments occurs and the covariance keeps its combinatorial interpretation, allowing us to enumerate a larger variety of graphs.

Finally, applying the same strategy than in the convex case, we are able to prove the convergence of the free energy.

THEOREM 7.6. *For L is bigger than some $L_0 > 0$, there exists $\eta$ such that, for $|\mathbf{t}| < \eta$, the following asymptotics hold:*

$$\log Z_{V_\mathbf{t}, L}^N = N^2 F^0(V_\mathbf{t}) + F^1(V_\mathbf{t}) + o(1),$$

*with*

$$F^0(V_\mathbf{t}) = \sum_{k_1, \ldots, k_n \in \mathbb{N} - \{0\}} \prod_{i=1}^n \frac{(-t_i)^{k_i}}{k_i!} \mathcal{M}_{k_1, \ldots, k_n}$$

*and*

$$F^1(V_\mathbf{t}) = \sum_{k_1, \ldots, k_n \in \mathbb{N} - \{0\}} \prod_{i=1}^n \frac{(-t_i)^{k_i}}{k_i!} \mathcal{M}_{k_1, \ldots, k_n}^1.$$

**Acknowledgments.** Part of this research was accomplished when Alice Guionnet was visiting the Miller institute for Basic Research in Science, University of California Berkeley and when Edouard Maurel–Segala was visiting the Department of Statistics at Stanford University.

ECOLE NORMALE SUPÉRIEURE DE LYON
UNITÉ DE MATHÉMATIQUES PURES ET APPLIQUÉES
UMR 5669 46 ALLÉE D'ITALIE
69364 LYON CEDEX 07
FRANCE
E-MAIL: aguionne@umpa.ens-lyon.fr
emaurel@umpa.ens-lyon.fr